\documentclass[aos]{imsart}

\RequirePackage{amsthm,amsmath,amsfonts,amssymb}
\RequirePackage[numbers]{natbib}
\RequirePackage[colorlinks,citecolor=blue,urlcolor=blue]{hyperref}
\RequirePackage{graphicx}

\usepackage{caption}
\usepackage{subcaption}
\usepackage{bm}

\startlocaldefs
\theoremstyle{plain}

\newtheorem{theorem}{Theorem}[section]
\newtheorem{lemma}[theorem]{Lemma}

\newcommand{\la}{\left\{}
\newcommand{\ra}{\right\}}
\newcommand{\mb}{\mathbb}
\newcommand{\ve}{\varepsilon}

\newcommand{\rar}{\rightarrow}

\newcommand{\lb}{\left[}
\newcommand{\rb}{\right]}
\newcommand{\ds}{\displaystyle}
\newcommand{\mc}{\mathcal}
\newcommand{\lr}{\left(}
\newcommand{\rr}{\right)}

\newtheorem{thmy}{Theorem}
\newtheorem{defn}{Definition}

\theoremstyle{remark}

\newtheorem{thm}{Theorem}

\endlocaldefs

\begin{document}

\begin{frontmatter}
\title{Asymptotic Distributions of Largest Pearson Correlation Coefficients under Dependent Structures}
\runtitle{Distribution of Pearson Correlation Coefficient}

\begin{aug}
\author[A]{\fnms{Tiefeng} \snm{Jiang}\ead[label=e1]{jiang0410@umn.edu}}
\and
\author[B]{\fnms{Tuan} \snm{Pham}\ead[label=e3]{pham0310@umn.edu}}
\address[A]{Department of Statistics, University of Minnesota, Twin Cities \printead{e1}}

\address[B]{Department of Statistics, University of Minnesota, Twin Cities \printead{e3}}
\end{aug}

\begin{abstract}
Given a random sample from a multivariate normal distribution whose covariance matrix is a Toeplitz matrix, we study the largest off-diagonal entry of the sample correlation matrix. Assuming the multivariate normal distribution has the covariance structure of an auto-regressive sequence, we establish a phase transition in the limiting distribution of the largest off-diagonal entry. We show that the limiting distributions are of Gumbel-type (with different parameters) depending on how large or small the parameter of the autoregressive sequence is. At the critical case, we obtain that the limiting distribution is the maximum of two independent random variables of Gumbel distributions. This phase transition establishes the exact threshold at which the auto-regressive covariance structure behaves differently than its counterpart with the covariance matrix equal to the identity. Assuming the covariance matrix is a general Toeplitz matrix, we obtain the limiting distribution of the largest entry under the ultra-high dimensional settings: it is a weighted sum of two independent random variables, one normal and the other following a Gumbel-type law. The counterpart of the non-Gaussian case is also discussed. As an application, we study a high-dimensional covariance testing problem.   
\end{abstract}

\begin{keyword}[class=MSC]
\kwd[Primary ]{62E20}
\kwd[; secondary ]{62H15}
\end{keyword}

\begin{keyword}
\kwd{high-dimensional testing}
\kwd{extreme values}
\kwd{sample correlation matrix}
\kwd{uniform central limit theorem}
\end{keyword}

\end{frontmatter}

\section{Introduction}
Let  $\textbf{x}_1, \textbf{x}_2, \cdots, \textbf{x}_n$ be a random sample from a  $p$-dimensional population distribution with mean  $\bm{\mu}$ and covariance matrix $\bm{\Sigma}$. Write $\textbf{x}_i=(x_{i1},x_{i2},\cdots,x_{ip})^T$ for each $i$. The Pearson correlation $\hat{r}_{ij}$ between the $i$-th and $j$-th columns of the data matrix $(\textbf{x}_1, \textbf{x}_2, \cdots, \textbf{x}_n)^T$ is given by 
\begin{eqnarray}\label{haorizi}
\hat{r}_{i j}=
\frac{\sum_{k=1}^{n}\left(x_{k i}-\bar{x}_{i}\right)\left(x_{k j}-\bar{x}_{j}\right)}{\sqrt{\sum_{k=1}^{n}\left(x_{k i}-\bar{x}_{i}\right)^{2}} \sqrt{\sum_{k=1}^{n}\left(x_{k j}-\bar{x}_{j}\right)^{2}}},
\end{eqnarray}
where $\bar{\textbf{x}}_i=n^{-1}\sum_{k=1}^{n} x_{ki}$. Then $\bm{R}:=(\hat{r}_{i j})$ is the $p\times p$ sample correlation matrix. Assuming the $p$ entries of the population are i.i.d., that is, $x_{11},x_{12},\cdots,x_{1p}$ are i.i.d., Jiang \cite{Jiang04} shows that $\max_{1\leq i < j \leq p}\hat{r}_{i j}$ asymptotically follows a Gumbel distribution:
\begin{eqnarray}\label{little}
n \big(\max_{1\leq i < j \leq p}\hat{r}_{i j}\big)^{2}-4 \log p+\log \log p \ \ \mbox{goes to a Gumbel distribution}
\end{eqnarray}
weakly with cdf $\exp (-(8 \pi)^{-1/2} e^{-t / 2})$, provided 
 $\lim p/n = \gamma>0$ and $\mb{E}(|x_{ij}|^{30+\ve})< \infty$. The follow-up works focus on both theories and applications. We elaborate these next.
 
 On the theoretical side,  \eqref{little} was extended by Zhou \cite{Zhou07} to a more general setting which only requires finiteness of the $6$th moment of $x_{ij}$'s and $\ds{p=O(n)}$. In subsequent papers  \cite{Li1, Li2, Li3, ShaoZhou}, sufficient and necessary conditions are  provided for \eqref{little} to hold.
Regarding the dependence between $p$ and $n$, it was shown in Liu et al. \cite{Liu16} and Shao and Zhou  (2014) that \eqref{little} still holds when $p$ grows in a polynomial rate of $n$ such that $p =o(n^\alpha)$ for some $\alpha>0$. Cai and Jiang \cite{Jiang12} consider a weakly dependent setting with ultra high-dimensional scenario: $\log p = o(n^\beta)$ for some $\beta>0$ and obtain \eqref{little}. Some other work
are the study of \eqref{little} under the assumption that population has a spherical distribution (Cai and Jiang \cite{JiangCai12}) or the asymptotic distribution of the maximum pairwise geodesic distances (Cai et al. \cite{Jiang13}). In these two papers, interesting  transition phenomena are found and the asymptotic distribution depends on $c=\lim_{n\to\infty}(\log p)/n$ with $c=0$, $c\in (0, \infty)$ and   $c=\infty$.

On the application side, test statistics based on maximum of sample correlation coefficients have  been proved to perform well under sparsity assumptions. One can see this from,  for example,  Cai et al. \cite{Cai17, Cai2011, Cai2013, Cai2014, Cai16} and Feng et al. \cite{FengJiang}, in which the test powers are higher than those of other tests. In fact,  Cai et al. \cite{Cai2013} justifies this in their framework. In a different context, Chen and Liu \cite{ChenLiu} studies independence testing problem with ultra high-dimensional correlated samples. Their proposed test statistic is an analog of the maximum sample correlations and involves a consistent estimator of the covariance matrix's Frobenius norm. Another interesting development is to understand the performance of maximum sample correlation coefficients under non-sparse alternatives. Yu et al. \cite{Yu1, Yu2} investigate power-boosting properties of the test statistic in  Cai et al. \cite{Cai2013} by combining it with a statistic of quadratic form. This new statistic has much better performance against dense alternatives than the original one. On the other hand, Cai and Ma \cite{CaiMa13} study optimal tests by using the U-statistics based tests;  Li and Xue \cite{D.Li, D.Li2} studies asymptotic independence of U-statistics based on sample covariance matrix and maximum sample correlation coefficients, as well as their applications to covariance testing.

Recently, Fan and Jiang \cite{Jiang19}   $\max_{1\leq i < j \leq p}\hat{r}_{i j}$ appeared in \eqref{little} with the population distribution being Gaussian and  the population correlation coefficient $r_{ij}=r$ for any $i, j$, where $r$ depends on $n$ only. They find there is a phase transition in the limiting distribution of $\max_{1\leq i < j \leq p}\hat{r}_{i j}$. The phase transition occurs as $r \sqrt{\log p} \rightarrow c$ with $c$ being finite. The limiting distributions according to $c=0$, $c=\infty$ and $c\in (0, \infty)$ are the Gumbel distribution, the normal distribution and the convolution of the two, respectively. The latter is also the the distribution of two independent random variables: one is normal and the other is Gumbel. 

In this paper, we will study the asymptotic distribution of $\max_{1\leq i < j \leq p}\hat{r}_{i j}$  as the population distribution is Gaussian and its covariance matrix $\bm{\Sigma}$ has two special structures, and then briefly discuss the non-Gaussian case and give applications. First, the two features of  $\bm{\Sigma}$ are given as follows. 

\begin{itemize}
      \item 
$\bm{\Sigma}$ is the covariance matrix of the auto-regressive model $AR(1)$, that is, 
     \begin{align} \label{cov1}
         \bm{\Sigma}=(r^{|i-j|})_{p\times p}=\left(\begin{array}{cccc} 
1 & r & \cdots & r^{p-1} \\
r & 1 & \cdots & r^{p-2} \\
\vdots & \vdots & & \vdots \\
r^{p-1} & r^{p-2} & \cdots & 1
\end{array}\right),~~~ 0\leq r <1;
     \end{align}
    
\item 
$\bm{\Sigma}$ is a Toeplitz matrix associated with a \textit{fixed}, non-increasing sequence $\la r_k \ra_{k=1}^{\infty}$, that is, 
  \begin{align} \label{cov2}
      \bm{\Sigma}=(r_{|i-j|})_{p\times p}=\left(\begin{array}{cccc}
1 & r_1 & \cdots & r_{p-1} \\
r_1 & 1 & \cdots & r_{p-2} \\
\vdots & \vdots & & \vdots \\
r_{p-1} & r_{p-2} & \cdots & 1
\end{array}\right)
  \end{align}
where $1=r_0 \geq r_1 \geq r_2 \geq \cdots \geq 0$. 
\end{itemize}
The matrix in \eqref{cov2} is the covariance matrix of a stationary sequence of random variables. Superficially,  \eqref{cov1} is a special case of \eqref{cov2}. The difference is that the quantity ``$r$" in \eqref{cov1} may change with $n$, but ``$r_1, r_2, \cdots$" in \eqref{cov2} are free of $n$ and $p$, which is the reason we put the term ``{\it fixed}". Indeed, if ``$r$" in \eqref{cov1} remains fixed, then \eqref{cov1} is a special case of  \eqref{cov2}.  

In the study, we will consider the ultra high-dimensional scenario, in which $\log p$ can grow as fast as $n^c$ for some constant $c>0$.

As the Pearson correlation coefficients are invariant under translation and scaling, we may assume, without loss of generality, that mean vector of the Gaussian population is $\textbf{0}$. Moreover, it is well-known that under the Gaussian assumption and the positive definiteness of $\bm{\Sigma}$ (see Fan and Jiang \cite{Jiang19}), the statistic $\max_{1\leq i < j \leq p}\hat{r}_{i j}$ has the same distribution as
$$
\max_{1\leq i<j \leq p} \frac{\sum_{k=1}^{n-1} x_{ki}x_{kj}}{\sqrt{\sum_{k=1}^{n-1} x_{ki}^2 \cdot \sum_{k=1}^{n-1} x_{kj}^2}}.$$
For the sake of simplicity, we replace ``$n-1$" above by ``$n$" and define
\begin{eqnarray}\label{wangshi}
\hat{\rho}_{i,j}=\frac{\sum_{k=1}^{n} x_{ki}x_{kj}}{\sqrt{\sum_{k=1}^{n} x_{ki}^2\sum_{k=1}^{n} x_{kj}^2}}.
\end{eqnarray}
We will work on the statistic
\begin{eqnarray}\label{language}
\mc{L}_n =\max_{1\leq i<j \leq p} \hat{\rho}_{i,j}.  
\end{eqnarray}

We will see $\mc{L}_n$ in \eqref{cov1} and \eqref{cov2} behave very differently. This is due to the assumption that $r$ from \eqref{cov1} depends on $n$ but $r_1, r_2, \cdots$ from \eqref{cov2} are free of $n$ and $p$. In fact $\mc{L}_n$ has phase transitions in both cases, however, their behaviors are distinct.  For case \eqref{cov1}, $\mc{L}_n$ is asymptotically Gumbel for all $r$ except for a critical value $r$, in which case the limit is the maximum of two independent Gumbel-distributed random variables. For case \eqref{cov2},  the limit of $\mc{L}_n$ is not Gumbel in big regimes.

\indent For the proofs of our results, the techniques  employed here is different from the Stein  Poisson approximation method used in the earlier works aforementioned. 
Instead, we compare the distribution of the maximum of Pearson correlation coefficients with the maximum of a Gaussian random field. This step is taken care of by the uniform CLT developed by Chernozhukov et al. \cite{Kato20}. Then the  problem is reduced to studying the maximum of a (nonstationary) Gaussian sequence.  

The rest of the paper is organized as follows. 
In Section \ref{sec3}, we state our  main results and contributions. We also discuss  the non-Gaussian case. In in Section \ref{sec4}, simulation and an application to high dimensional test are provided. Section \ref{discussion} contains some concluding remarks. The proofs of the main results are given in Section \ref{sec6}.

\section{Main results and discussion} \label{sec3}
In this section, we will present our main results in the order of notation, the behavior of $\mc{L}_n$ from \eqref{language} under the $AR(1)$ and Toeplitz structures, respectively. We always assume  $\textbf{x}_1, \textbf{x}_2, \cdots, \textbf{x}_n$ are i.i.d. $\mathbb{R}^p$-valued random variables with distribution $N(\bm{0}, \bm{\Sigma})$.
Finally we will make a brief discussion on  the non-Gaussian case.

\subsection{Notation} \label{sec2}
Throughout the paper, we assume that $p\geq 3$ and $n\geq 3$. For a sequence of vectors $\bm{x}_k \in \mb{R}^{p}$, we will write $\bm{x}_k=(x_{k1}, \cdots, x_{kp})^T$ or  $\bm{x}_k=(x_{k,1}, \cdots, x_{k,p})^T$ for clarity. We sometimes also write $x_k(i)$ for  $x_{ki}$ or $x_{k,i}$. 

Given random variable $X$, let $\|X\|_{q}=(\mb{E}|X|^q)^{1/q}$ be the standard $L_{q}$ norm for $q\geq 1$. We also define the Orlicz norm $||.||_{\psi_q}$ by 
\begin{eqnarray}\label{haoshi}
||X||_{\psi_q}=\inf \big\{t>0: \mb{E}\exp\big(|X|^q/t^q\big) \leq 2 \big\}.
\end{eqnarray}

Let $I\subset \{1, 2, \cdots\}^2$ be an index set. With a slight abuse of notation, we will sometimes write $X_{i,j}$ or $X(i,j)$ to indicate the $(i,j)$-coordinate for random vector $X=(X_{i,j})_{(i,j) \in I}$. \\
\indent Unless stated otherwise, the quantities $r,p$ all depend on $n$ and $p=p_n \rightarrow \infty$.\\
\indent For vectors $x=(x_1, \cdots, x_d)^T\in \mathbb{R}^d$ and $y=(y_1, \cdots, y_d)^T\in \mb{R}^d$, we use the notation $x \leq y$ to indicate  $x_i \leq y_i$ for each  $1\leq i \leq d$. 

For a random sequence $\xi_n$ and a nonrandom sequence $a_n, \xi_n=o_{\mb{P}}\left(a_n\right)$ means $\xi_n / a_n \rightarrow 0$ in probability as $n \rightarrow \infty$; $\xi_n=O_{\mb{P}}\left(a_n\right)$ means $\lim _{C \rightarrow \infty} \limsup _{n \rightarrow \infty} P\left(\left|\xi_n / a_n\right|>C\right)=0$.

\subsection{The $AR(1)$ covariance structure}
In this part we will work on $\mc{L}_n$ in  \eqref{language} under the $AR(1)$ structure. If there is no confusion we will write $r=r_n$ for convenience. 

\begin{thm}\label{slowdecay}
Assume $\log p =o(n^{1/5})$ and $\limsup_{n\to\infty} r_n < 1$. Assume also

$$\lim_{n \rightarrow \infty} \frac{r \sqrt{n}}{\sqrt{\log p}}=L\in [0, \infty].$$
The following hold. 
(i). If $0 \leq L < 2 -\sqrt{2}$, then  $a_n \mc{L}_n - b_n$
converges weakly to a Gumbel distribution with cdf $F(x)=e^{-K_1 e^{-x}}$, where 
\begin{eqnarray*}
a_n = 2\sqrt{n \log p},\  \ b_n  = 4 \log p - \frac{1}{2} \big[ \log \log p + \log(4 \pi) \big] \ \ \mbox{and}\ \ K_1 &= \frac{1}{2 \sqrt{2}}.
\end{eqnarray*}

(ii).    If $ 2-\sqrt{2} < L \leq \infty $, then 
$c_n \mc{L}_n -d_n$
converges weakly to a Gumbel distribution with cdf $F(x)=e^{- e^{-x}}$, where
\begin{eqnarray*}
c_n = \frac{\sqrt{2n \log p}}{1-r^2}\ \ \ \mbox{and}\ \ \ d_n = \frac{r \sqrt{2n \log p}}{1-r^2} + 2\log p - \frac{1}{2} \big[ \log \log p + \log(4 \pi) \big]. 
\end{eqnarray*}
\end{thm}

Recalling \eqref{language}, $\mc{L}_n$ is the maximum of $p(p-1)/2$ random variables. The above implies that, if $r$ decays fast enough, that is,  $L<2-\sqrt{2}$, then $\mc{L}_n$ behaves like the maximum of $p(p-1)/2$ i.i.d. standard normals. However, as $L>2-\sqrt{2}$, the maximum $\mc{L}_n$ behaves like the maximum of $p$ i.i.d.  normal random variables. Theorem \ref{slowdecay}(i) includes the case $r=0$, that is, $\bm{\Sigma}=\bm{I}_p$. Obviously, this says that there is a phase transition at $L=2-\sqrt{2}.$

 The above transition phenomenon can be explained in the following way. Recall $\mc{L}_n =\max_{1\leq i<j \leq p} \hat{\rho}_{i,j}$. If $r=r_n$ is smaller than the threshold, then  every $\hat{\rho}_{i,j}$ with  $1\leq i<j \leq p$ contributes to $\mc{L}_n$ equally. However, if $r=r_n$ is larger than the threshold, only $\hat{\rho}_{i,i+1}$ with $1\leq i \leq p-1$ essentially contribute to $\mc{L}_n$. This is caused by the fact $\hat{\rho}_{i,j} \sim r^{|j-i|}$ by the law of large numbers. As a consequence,  $\hat{\rho}_{i,i+1} \sim r$ and the rest $\hat{\rho}_{i,j}$ are of orders $r^2, r^3, \cdots, r^{p-1}$. Hence, if the difference between $r$ and $r^2$ is not too small, the maximum $\mc{L}_n$ is achieved at $\hat{\rho}_{i,i+1}$ with $1\leq i \leq p-1$   with high probability.

We prove Theorem \ref{slowdecay} by comparing the distribution of the maximum of $\hat{\rho}_{i,j}$ with the maximum of a Gaussian random field. The argument is carried out by the uniform CLT derived by Chernozhukov et al. \cite{Kato20}. Then we bring the problem down to the study of the maximum of a nonstationary Gaussian sequence. While we believe Theorem \ref{slowdecay} can still be proven by the Stein method for Poisson approximation, employed in the earlier works aforementioned, a preliminary analysis shows 
that the argument is clumsy and lengthy. The latter is due to the computation of conditional probabilities.

Recall $L$ in Theorem \ref{slowdecay}. Our next result investigates $\mc{L}_n$ at the critical case $L=2-\sqrt{2}$. A finer analysis yields the following result.
\begin{thm} \label{large}
Assume the setting in Theorem \ref{slowdecay} with $L=2-\sqrt{2}$. Then
$$
{ \kappa_n:}= \frac{r \sqrt{n}}{\sqrt{\log p}} -(2-\sqrt{2}) \to 0.
$$
Define 
$
\lambda_n=\sqrt{2} (\log p) \kappa_n + (8^{-1/2}-2^{-1})\log \log p.
$
Let $K_1$, $a_n$, $b_n$, $c_n$ and $d_n$ be defined as in Theorem \ref{slowdecay}. The following statements hold.

(i). If $\lambda_n \rar -\infty$ then 
$a_n \mc{L}_n - b_n$
converges weakly to a Gumbel distribution with cdf $F(x)=e^{-K_1 e^{-x}}.$ 

(ii).  If $\lambda_n \rar \lambda$ then 
    $a_n \mc{L}_n - b_n$
converges weakly to a probability distribution with cdf 
$$
 F(x)=\exp\big(-K_1e^{-x}\big) \cdot \exp\big(-e^{-(x/\sqrt{2})-K_2+\lambda}\big),
 $$
    where $K_2 = (\frac{1}{2} - \frac{1}{\sqrt{8}})\log (4 \pi)$. Obviously, $F(x)$ is the distribution of the maximum of two independent random variables with different Gumbel-type distributions. 

(iii). If $\lambda_n \rar \infty$ then $c_n \mc{L}_n - d_n$ 
converges weakly to a Gumbel distribution with cdf 
$F(x)= e^{- e^{-x}}$.
\end{thm}
\indent We have seen a subtlety about the limiting distribution of $\mc{L}_n$. It changes according to how ${\kappa _n}$ goes to zero. 
The phase transition phenomenon occurred in Theorems \ref{slowdecay} and \ref{large} is interesting. We now make a quick summary as follow. 
\begin{itemize}
    \item When the ratio $\frac{r \sqrt{n}}{\sqrt{\log p}}$ is above $2-\sqrt{2}$, only the first sub-diagonal entries of the sample correlation matrix $(\hat{\rho}_{i,j})$ contribute to the maximum $\mc{L}_n$.
    \item When the ratio $\frac{r \sqrt{n}}{\sqrt{\log p}}$ is below $2-\sqrt{2}$, in contrast to the previous case, the contribution of the first sub-diagonal entries of  $(\hat{\rho}_{i,j})$ is negligible and the rest entries make the whole contribution to $\mc{L}_n$.

    \item In the critical case for which the ratio $\frac{r \sqrt{n}}{\sqrt{\log p}}$ is approximately $2-\sqrt{2}$, a careful analysis shows that the limiting distribution of $\mc{L}_n$ could fall into either of the two previous cases or it can be the maximum of two independent Gumbel distributions.
\end{itemize}

\medskip

The following result studies $\mc{L}_n$ under the case that $r=r_n\to 1$.

\begin{thm} \label{fastdecay}
Assume the setting in Theorem \ref{slowdecay} with $\log p = o(n^{1/7})$. Recall the notation $c_n$ and $d_n$. If $r=r_n\to 1$ and $1-r \geq C/\log p$ for a positive constant $C$ free of $n$, then $c_n\mc{L}_n-d_n$ converges weakly to a Gumbel law with cdf $\exp\{-e^{-x}\}$.
\end{thm}

The proof of Theorem \ref{fastdecay} is provided in the supplement Jiang and Pham \cite{JiangPham21}.

The above result asserts that when $r$ converges to $1$ at a rate not faster than  $1/\log p$, similar to conclusion (i) of Theorem \ref{slowdecay}, the limiting distribution of $\mc{L}_n$ is still a Gumbel distribution. Notice Theorem \ref{slowdecay} holds under restriction $\limsup_{n\to\infty} r_n < 1$. If $\limsup_{n\to\infty} r_n = 1$, Theorem \ref{fastdecay} says that  the limiting distribution of $\mc{L}_n$ still exists and it is a Gumbel law as long as $r=r_n\to 1$ not so fast. It is possible that the same result holds as $r\to 1$ at a faster rate. However, we will not pursue such a technical improvement in this paper.

On the other hand, we could combine Theorem \ref{fastdecay} with  (ii) of Theorem \ref{slowdecay}
 with $L=\infty$ in a single conclusion if we ignore the assumptions $\log p = o(n^{1/7})$ in the former theorem and $\log p = o(n^{1/5})$ in the later theorem. It is conceivable that a more stringent condition is needed for the extreme case that $r=r_n\to 1$.

\subsection{The Toeplitz covariance structure}
Let $1=r_0 > r_1 \geq r_2 \geq\cdots$ be a fixed sequence of  non-negative numbers and $r_i>r_j$ for some $j>i\geq 1$. To have a meaningful problem, we assume the Toeplitz covariance matrix $\bm{\Sigma}_p$ formed by $\la 1,r_1,r_2,\cdots,r_{p-1} \ra$ given in \eqref{cov2} is positive definite for all $n\geq 1$. A sufficient condition to guarantee that $\bm{\Sigma}_n$ is positive definite is provided by Polya's criterion in \cite{Polya}.
Let $d \geq 1$ be the smallest integer such that $r_d > r_{d+1}$. We now study $\mc{L}_n$ in the following.
\begin{thm} \label{Toeplitz}
Assume $\log p = o(n^{1/5})$.  Then the following hold.

(i) If $r_p \sqrt{\log p} \rightarrow \gamma \in [0, \infty)$ then $ a_n^{*}\mc{L}_n-b_n^{*}$ converges weakly to  $-\gamma_0 + \sqrt{2\gamma_0} Z + G-(1/2)\log (4\pi)$, where $\gamma_0=2\gamma^2(1+r_1)^{-2}, 
~~Z \sim N(0,1),~~ G \ \mbox{has}\ \mbox{cdf}\  \exp(-{e^{-x}}),$  $Z$ and $G$ are independent, and  
\begin{eqnarray*}
a_n^{*}&=&\frac{1}{1-r_1^2}\sqrt{2n\log (pd)};\,\\
b_n^{*}& =&\frac{r_1}{1-r_1^2}\sqrt{2n\log (pd)}+ 2 \log (pd)-\frac{1}{2}\log \log (pd).
\end{eqnarray*}

(ii) Assume that $r_p \to 0$, $r_p \sqrt{\log p} \to \infty$ and $r_{k+1}/r_k \to 1$. Define $f(0)=1$ and
$$f(k)= \frac{1}{(1-r_1^2)^2} \Big[ r_1^2 r_k^2 + \frac{1}{2}r_1^2(r_{k-1}^2 + r_{k+1}^2) + r_k^2 + r_{k-1}r_{k+1} -2r_1r_k(r_{k-1}+r_{k+1}) \Big]$$
for all $k\geq 1$. Assume $d=1$, $f(k) \log k$ increases  to infinity and there exists  $K\geq 1$ such that 
$f(k)$ is non-increasing in $[K, \infty)$. Then, as $n\to\infty$, 
$$ \frac{\sqrt{n}\mc{L}_n -r_1 \sqrt{n}}{(1-r_1^2)\sqrt{f(p-1)}}  - \sqrt{\frac{1-f(p-1)}{f(p-1)}} \cdot \left(\sqrt{2 \log p} - \frac{\log \log p+ \log(4 \pi)}{2\sqrt{2 \log p}} \right) \to N(0, 1).$$

\end{thm}

Let us briefly explain the role of index $d$ in the statement of Theorem \ref{Toeplitz}. Recall $\mc{L}_n =\max_{1\leq i<j \leq p} \hat{\rho}_{i,j}$. 
By the law of large numbers, $\hat{\rho}_{ij}$ concentrates at $r_{|i-j|}$. We find that $\mc{L}_n$  is attained in the subset ${\la (i,j):\, 1\leq i<j,\, j-i \leq d \ra}$ with a high probability. Similar to the discussion below Theorem \ref{slowdecay}, this is due to the gap between $r_d$ and $r_{d+1}, r_{d+2}, \cdots.$  Therefore, $\mc{L}_n$ is roughly the maximum of a Gaussian array of size $(p-1) +(p-2)+\cdots +(p-d) \approx pd$. This observation also explains why the normalizing constants $a_n^*$ and $b_n^*$ have such forms.

In the statement of Theorem \ref{Toeplitz}(ii), the function  $f(k)$ is required to satisfy certain conditions.  Example 5 below shows that a common Toeplitz matrix satisfies those restrictions. Although  function $f(k)$ does not show a friendly look to reveal its decreasing  property, we can use computing softwares, for example, MATLAB to finish the job easily.  We only need to check the dominated term of $f(k)-f(k+1)$  is non-negative. Examples and codes for MATLAB are provided in  Section 3.5 from the supplement \cite{JiangPham21}. 

We assume the sequence $\la r_n \ra_{n \geq 1}$ is a fixed sequence in  Theorem \ref{Toeplitz}. It is possible to prove a slightly stronger result by letting the sequence $\{r_1, \cdots, r_n\}$ depend on $n$ for each $n\geq 1$. However, one must impose extra conditions on the regularity of these sequences and also the growth of $d$ to ensure similar results.  
 In the special case when $\gamma=0$, the asymptotic distribution is the classical Gumbel distribution which have been obtained in similar settings (see Jiang \cite{Jiang04}, Zhou \cite{Zhou07}, Li et al. \cite{Li1,Li2,Li3}, Cai and Jiang \cite{Jiang12}). However, in the case $\gamma>0$, the dependence between the numerators and denominator in \eqref{wangshi} is strong and it is no longer negligible. The technique we used here can not handle the case $\gamma=\infty$ and $d \geq 2$ due to the non-stationarity and the high dependence among $\{\hat{\rho}_{i,j}, j>i\geq 1\}$. We hope new techniques are created to tackle the case.

 Item (ii) of Theorem \ref{Toeplitz} is an analog of a classical result in the extreme value theory regarding for strongly dependent, stationary Gaussian sequence (see Theorem 6.6.4 in \cite{Leadbetter}). Let us discuss the assumptions imposed in the theorem. The condition $d=1$ is needed to transform the problem to a stationary scenario. The conditions $r_{k+1}/r_k \to 1$ and $r_k$ being non-increasing are required to eliminate the possibility that $\la r_k \ra_{k \geq 1}$ behaves irregularly. The regularity behaviour of $\la r_k \ra_{k\geq 1}$ should be expected as it is the indeed the case in most statistical applications. We would like to point out that the irregular behaviour of $\la r_k \ra_{k\geq 1}$ can lead to non-standard asymptotic distributions. For example, generate two independent sequences  $\la z_1, z_3, z_5,... \ra$ and $\la z_2, z_4, z_6,... \ra$ such that each one is a Gaussian sequence with equicorrelated correlation $1/3$. Merge the two sequences to  form  $\{z_1,z_2, z_3, z_4, ...\}$. Under an appropriate scaling, we can see that $M_n:= \max_{1 \leq i \leq n} z_i$  is asymptotically the maximum of two i.i.d. $N(0, 1)$'s. As in the proof of Theorem \ref{Toeplitz}, this phenomenon also holds for $\max_{1 \leq i \leq p-1} \hat{\rho}_{i,i+1}$, which is essentially identical to $\mc{L}_n$. Note that this distribution is no longer an extreme value distribution of type I, II or III.

\indent Let us finish this section by providing some examples on Toeplitz matrix $\bm{\Sigma}$ of the form (\ref{cov2}) and their connections to existing results in literature. 

\begin{itemize}
    \item {\it Example 1}. Consider the covariance matrices of stationary and $m$-dependent Gaussian sequences. Write $\bm{\Sigma}_{ij}=a_{|i-j|}$. Then  $a_n$ vanishes for all $n \geq m$.
    Such matrices are very common, practical covariance matrices of the form (\ref{cov2}). They also belong to the class of banded covariance matrices, which has received a lot of attention in the last ten years; see, for example, \cite{Jiang12, Jiang2019}. 
    In this case, Theorem \ref{Toeplitz} holds with $d=m_1$ (the smallest numer such that $a_{m_1}>a_{m_1+1}$),  $r_1=a_1$ and $\gamma=0$. 
    
    \medskip
    \item  {\it Example 2}. Consider the covariance matrices $\bm{\Sigma}$ of the form 
     \[ \bm{\Sigma}_{ij} = \begin{cases} 
          1, & i=j; \\
          \frac{a}{\sqrt{\log(A+|i-j|)}}, & |i-j| \geq 1
       \end{cases}
    \]
    where $a >0$ and $A \geq 1$. One can see that $\bm{\Sigma}$ is positive semi-definite by Polya's  criterion (see \cite{Polya}). Moreover, for $r_n= a \cdot (\log(A+n))^{-1/2}$, it holds that
    $\lim_{n \rar \infty} r_n \sqrt{\log n} = a. $
    Therefore, the conditions of Theorem \ref{Toeplitz}(i) are satisfied with $d=1$,  $r_1=a(\log (A+1))^{-1/2}$ and $\gamma=a$. The asymptotic distribution of $\mc{L}_n$ is the mixture of two independent distributions: the Gumbel and the standard normal.
    
    \medskip
    \item  {\it Example 3}. Consider the $AR(1)$ covariance matrix with fixed parameter $r$. In this case,   $d=1$, $r_1=r$ and $\gamma=0$. The limit  of $\mc{L}_n$ in Theorem \ref{Toeplitz} agrees with that of Theorem \ref{slowdecay}.

    \medskip
    \item {\it Example 4}.  Consider the $p\times p$ covariance matrix  $\bm{\Sigma}=\bm{\Sigma}_n$ from \eqref{cov2} with 
      \[ \bm{\Sigma}_{ij} = \begin{cases} 
          1, & i=j; \\
          r_1, & |i-j| = 1;\\
          r_{n}, & |i-j| \geq 2
       \end{cases}
    \]
    with $\max\{2r_1-1, 0\} \leq  r_n \leq r_1$ for each $n\geq 1$. Of course, this is equivalent to $0\leq r_n \leq r_1$ if $r_1\leq 1/2$. In other words, $\bm{\Sigma}_{ij}=r_1$ for $|i-j|=1$ and 
 $r_1$ does not depend on $n$; if $|i-j|\geq 2$, $\bm{\Sigma}_{ij}=r_n$ but not depend $i,j$.  It is checked in Section 3.3  from the supplement \cite{JiangPham21}  that $\bm{\Sigma}_n$ is positive definite. We assume $r_n \sqrt{\log n} \to \infty$ ($r_n$ is not required to go to $0$).  Note that this is not the exact setting presented in Theorem \ref{Toeplitz} since $\bm{\Sigma}$ depends on $n$. It can be checked that Theorem \ref{Toeplitz}(ii) still holds  with ``$f(p-2)$" replaced with ``$2(1+r_1)^{-2}r_{n}^2$". The proof remains almost the same with minor changes. We omit details. This is a variant of the results from Fan and Jiang \cite{Jiang19}. 
   \medskip
   \item {\it Example 5}.   Consider the class of covariance matrices $\Sigma$ satisfying the conditions of Theorem \ref{Toeplitz}(ii) with $d=1$ and $r_k =  [\log (k+2)]^{-1/2 + \ve}$ for $0< \ve <1/2$. It is verified in Section 3.4  from the supplement \cite{JiangPham21} that the corresponding sequence $f(k)$ decreases to $0$ and $f(k) \log k$ increases to $\infty$.   
    
\end{itemize}

\subsection{Potential extensions to non-Gaussian cases}
In this section we shall discuss  Theorems  \ref{slowdecay}, \ref{large}, \ref{fastdecay} and  \ref{Toeplitz} in the absence of Gaussian assumptions. It turns out that the results still hold, provided one has a $4$-th order moment matching condition on the distribution of the entries. We shall make it precise below.

\begin{defn} \label{def1}
For two random vectors (or two distributions) $\bm{X}=(x_1,x_2,...,x_n)$ and $\bm{Y}=(y_1,y_2,...,y_n)$, we say that $\bm{X}$ and $\bm{Y}$ satisfy the $4$-th order moment matching condition if the following equality holds for all $1\leq m \leq 4$ and $1 \leq i_1,i_2,...,i_m \leq n$:
$$\mb{E} \lr  \prod_{k=1}^{m} x_{i_k} \rr = \mb{E} \lr  \prod_{k=1}^{m} y_{i_k} \rr.$$
\end{defn}

Before stating the results, we make some assumptions on  the random sample $\bm{x}_1,\bm{x}_2,...,\bm{x}_n $.
   
    {\it Assumption 1}. $\bm{x}_i$'s are i.i.d. $p$-dimensional random vectors with sub-Gaussian tails: 
    $$\sup_{1 \leq k \leq p} \| x_{1k} \|_{\psi_2} \leq C,$$
    where $\bm{x}_i=(x_{i1},x_{i2},...,x_{ip})^T$, $\|.\|_{\psi_2}$ is the Orlicz norm in (\ref{haoshi}) and $C$ is an absolute constant.
    
     {\it Assumption 2.} For all $1 \leq k \leq p$, one has $\mb{E}x_{1k}=0$ and $\mb{E}x_{1k}^2=1$.

\medskip

Recall $\hat{r}_{i j}$, $\hat{\rho}_{i,j}$ and $\mc{L}_n$ in  \eqref{haorizi}, \eqref{wangshi} and  \eqref{language}, respectively.  We mention earlier that the maximum correlation coefficient $\tilde{\mc{L}}_n:=\max_{1\leq i < j \leq p}\hat{r}_{i j}$, under  normality assumptions, has the same distribution as $\mc{L}_{n-1}$. However, this may not be true in general. Even so, we still believe  the following result on non-Gaussian cases should hold (part of the verification is provided in \cite{JiangPham21}).

{\it Claim A}. 
{\it Suppose $\bm{x}_1,\bm{x}_2,...,\bm{x}_n$ are i.i.d. $p$-dimensional random vectors and that the entries of $\bm{x}_1$ satisfy Assumptions 1 and 2. Assume additionally that $\bm{x}_1$ and the multivariate normal $N_p(0,\bm{\Sigma}^{*})$ satisfy the moment matching condition described in Definition \ref{def1}.Then,
\begin{enumerate}
\item If $\bm{\Sigma}^{*}$ is an $AR(1)$ covariance matrix with corresponding parameter $r$ then Theorems  \ref{slowdecay}, \ref{large} and \ref{fastdecay} still hold for $\mc{L}_n$  and  $\tilde{\mc{L}}_n$, respectively;

\item If $\bm{\Sigma}^{*}$ is a positive semidefinite Toeplitz matrix formed by a positive, non-decreasing sequence $\la r_n \ra_{n \geq 1}$ then Theorems  \ref{Toeplitz} still holds for $\mc{L}_n$ and $\tilde{\mc{L}}_n$.
\end{enumerate}
}

\section{Applications and simulation results} \label{sec4}
\subsection{Two high-dimensional testing problems}
In this part, we consider an application of our main results to two problems of high-dimensional  testing for covariance matrices. Suppose $\bm{x}_1, \bm{x}_2,\cdots, \bm{x}_n$ is a random sample from the $p$-dimensional Gaussian population $N(\bm{\mu}, \bm{\Sigma})$ with known mean vector $\bm{\mu}$ and unknown covariance matrix $\bm{\Sigma}$.

\medskip 

\indent {\it Application 1: an independence test based on high dimensional data} 
    
    Consider the test with
    $H_0: \bm{\Sigma}=\bm{I}_p \ \ \ \mbox{vs}\ \ \ H_{1}:\bm{\Sigma} \neq \bm{I}_p.$
As mentioned in the introduction, for nonparametric testing problem involving with sparse alternatives, test statistics of extreme-value types tend to perform well. We will use the largest Pearson correlation coefficient as our test statistic:
$$\tilde{\mc{L}}_n = \max_{1 \leq i<j \leq p}\frac{\sum_{k=1}^{n}\left(x_{k i}-\bar{x}_{i}\right)\left(x_{k j}-\bar{x}_{j}\right)}{\sqrt{\sum_{k=1}^{n}\left(x_{k i}-\bar{x}_{i}\right)^{2}} \sqrt{\sum_{k=1}^{n}\left(x_{k j}-\bar{x}_{j}\right)^{2}}}.$$

\noindent As explained in \eqref{wangshi}, we know the test statistics $\mc{L}_{n-1}$ from \eqref{language} and  $\tilde{\mc{L}}_n$  have the same distribution. By taking $r_n=0$ for each $n\geq 1$,  it follows from (i) of Theorem \ref{slowdecay} that
$$2\sqrt{\log p}\Big(\sqrt{n -1}\tilde{\mc{L}}_{n} -2\sqrt{\log p} + \frac{\log \log p + \log 4\pi}{4 \sqrt{\log p}}\Big) \rightarrow G_1,$$
where $G_1$ is a Gumbel-type distribution with distribution function $F_{G_1}(x)= e^{-Ke^{-x}}$ and $K=\frac{1}{2\sqrt{2}}$. Reject $H_0$ if $\tilde{\mc{L}}_{n}$ is large. For a given size $\alpha>0$, let $q_{\alpha}$ be given by
$$q_{\alpha}=\log K -\log \log \frac{1}{1-\alpha}. $$

\noindent With this choice of $q_{\alpha}$, one has $\mb{P}(G_1 > q_{\alpha})=\alpha$. Let the critical value $c_n(\alpha)$ be chosen as 
$$c_n(\alpha)= \frac{q_{\alpha}}{2\sqrt{ \log p}} + 2\sqrt{ \log p} - \frac{\log \log p + \log 4\pi}{4\sqrt{\log p}}.$$

\noindent We then reject $H_0$ if $\sqrt{n{ -1}}\tilde{\mc{L}}_{n} > c_n(\alpha)$. We now analyze the powers of the proposed test under two  alternatives: one based on $AR(1)$ with covariance matrix given at \eqref{cov1}, and one based on the Toeplitz covariance matrix given at \eqref{cov2}.

Consider a  $AR(1)$ alternative $\bm{\Sigma}_n$ with associated parameters $r=r_n$ and assume 
 $$\liminf \frac{r \sqrt{n-1}}{\sqrt{\log p}} > 2-\sqrt{2}.$$
This condition is easily satisfied for many choices of $r_n$. For example, one could take $r_n$ to be a positive constant less than $1$ or $r_n=n^{-c}$ with $0<c<2/5$. Thanks to Theorem \ref{slowdecay}, it holds that 
\begin{eqnarray}\label{jingguo}
\mb{P} \left( \sqrt{2 \log p} \left( \frac{\sqrt{n-1}( \mc{L}_n -r)}{1-r^2}- \sqrt{2 \log p}+\frac{\log \log p+ \log 4\pi}{2\sqrt{2 \log p}} \right) \leq x \right) \rightarrow e^{-e^{-x}}.
\end{eqnarray}

\noindent Let $\beta(r)$ be the power function, one has the following identity
\begin{align*}
    \beta(r_n) &= \mb{P}(\sqrt{n-1}\mc{L}_n  > c_n(\alpha) | \bm{\Sigma}_n)\\
    &= \mb{P} \left(   \sqrt{2 \log p} \left( \frac{\sqrt{n-1} (\mc{L}_n -r)}{1-r^2}- \sqrt{2 \log p}+\frac{\log \log p+ \log 4\pi}{2\sqrt{2 \log p}} \right) >h_{n}(\alpha) \right),
\end{align*}
where $h_n(\alpha)$ is defined as
\begin{eqnarray}\label{cell_phone}
h_n(\alpha)= 2\sqrt{\log p} \lb \frac{c_n(\alpha) - r\sqrt{n-1}}{1-r^2}  -\sqrt{2 \log p} + \frac{\log \log p + \log 4 \pi}{2\sqrt{2 \log p}} \rb.
\end{eqnarray}
As $\log p = o(n^{1/5})$ and $p \rightarrow \infty$, it is easy to see that $h_n(\alpha)\to -\infty$. It follows from  \eqref{jingguo}  that the power function $\beta(r_n) \rightarrow 1$ as $n \rightarrow \infty$.

\indent{\it Application 2: Testing for auto-regressive covariance structure} 

Consider the testing problem  
\begin{align*}
    H_0: \bm{\Sigma}=\bm{\Sigma}_n~ in ~(\ref{cov1})~~ \textit{for some $r \in [0, 1)$} ~~ \mbox{vs} ~~ H_1: \bm{\Sigma}\not=\bm{\Sigma}_n~ in ~(\ref{cov1})~.
\end{align*}
Let $\bm{x}_1,\bm{x}_2, ...,\bm{x}_n$ be a sequence of $p$-dimensional normal random variables with mean  $\bm{\mu}$ and covariance matrix $\bm{\Sigma}$. The sample mean vector is 
\begin{align*}
\bar{x} &= \frac{\bm{x}_1+\bm{x}_2+...+\bm{x}_n}{n}=(\bar{x}_1,\bar{x}_2,...,\bar{x}_p)^T.
\end{align*}
Recall the notation at the beginning of Section \ref{sec2}. The $i$-th coordinate of the population distribution is estimated by the unbiased estimator
\begin{align*}
\hat{h}_i &= \frac{1}{n-1} \sum_{k=1}^{n} \lb  x_{k}(i) - \bar{x}_i \rb \cdot \lb x_{k}(i+1) - \bar{x}_{i+1} \rb.
\end{align*}
Now we estimate $r$ by using the average of sample covariances   over the first sub-diagonal: 
$$\hat{r} = \frac{1}{p-1} \sum_{i=1}^{p-1} \hat{h}_i. $$

It is checked in Section 3.2 from the supplement Jiang and Pham \cite{JiangPham21} that  $\mb{E}(\hat{r})=r$, $\hat{r}-r = O_{\mathbb{P}}(\frac{1}{\sqrt{np}})$, and  under $H_0$, 
\begin{eqnarray*}
~~\sqrt{2\log p} \lr \frac{ \sqrt{n} \mc{L}_n  -\hat{r}\sqrt{n}}{1-\hat{r}^2} - \sqrt{2\log p}+ \frac{\log \log p + \log 4 \pi}{2 \sqrt{2\log p}} \rr ~~\mbox{converges to} ~~ F(x)
\end{eqnarray*}
in distribution, where $F(x)$ is the Gumbel distribution with cdf $e^{- e^{-x}}.$  Based on this, 

we reject $H_0$ when $\sqrt{n} \mc{L}_n > s_n(\alpha)$, where the critical value $s_n(\alpha)$ is defined by 
$$s_n(\alpha)= \hat{r}\sqrt{n} + (1-\hat{r}^2) \cdot \left( \frac{q^{*}_{\alpha}}{\sqrt{2 \log p}} + \sqrt{2 \log p} - \frac{\log \log p + \log 4\pi}{2\sqrt{2\log p}} \right),$$
where $q^{*}_{\alpha}$ is the $(1-\alpha)$ quantile of the Gumbel distribution with CDF $F(x)=e^{-e^{-x}}$. 
Evidently, 
$\mb{P} \left( \sqrt{n} \mc{L}_n > s_n(\alpha) \right) \rightarrow \alpha.$
Therefore, the test is asymptotically of size $\alpha$.  Now we study the power under $H_a: \bm{\Sigma}=\bm{\Sigma}_n$, where $\bm{\Sigma}_n$ is as in \eqref{cov2}. We assume $r_{p} \sqrt{\log p} \rightarrow \gamma>0$. In this case, $r_1$ decides the limiting distribution completely. Similar to the computation in \eqref{cell_phone}, we obtain that the power function tends to $1$.

\indent{\it Application 3: Extreme angles of dependent random points on high-dimensional spheres.}

 Understanding geodesic distance of random points drawn uniformly on the hypersphere is an important problem in directional statistics, in which the direction of the data is of interests rather than its magnitude. In particular, the extreme geodesic distance has received a lot of attention in the last decade (see, for example, \cite{Jiang13} and the references therein). To the best of our knowledge, there has been no result concerning dependent random points in literature. Our main results shed light on  the behaviour of the extreme angles. To be more precise, suppose  $\bm{Y_1}, \bm{Y_2},...,\bm{Y_p}$ are drawn (not necessary independently) uniformly from  the sphere $\mathbb{S}^{n-1}$. It is well-known that 
\[
(\bm{Y_1}, \bm{Y_2},...,\bm{Y}_p) \stackrel{d}{=} \Big( \frac{\bm{X}_1}{\| \bm{X}_1 \|}, \frac{\bm{X}_2}{\| \bm{X}_2 \|},...,\frac{\bm{X}_p}{\| \bm{X}_p \|} \Big)
\]
where $X_i \sim N( \bm{0}, \bm{I}_n)$ for all $1 \leq i \leq p$. Note that $\bm{X}_i$'s are not independent. Then
\[ 
\max_{1 \leq i<j \leq p} \cos(\bm{Y_i}, \bm{Y_j}) \stackrel{d}{=} \mc{L}_n
\]
where $\mc{L}_n$ is defined in (\ref{language}) and $\bm{\Sigma}$ is  the covariance matrix of $(X_{11},..., X_{p1})^T$. We assume  $\bm{\Sigma}$ takes the form  (\ref{cov1}) or (\ref{cov2}). From Theorems \ref{slowdecay}-\ref{Toeplitz}, it is easy to determine the asymptotic distribution of the largest angle between the points $\bm{Y_i}$'s, which is the inverse cosine of a Gumbel, maximum of two independent Gumbels, mixture of Gumbel and normal or a standard normal depending on the form of $\bm{\Sigma}$ and the corresponding assumptions in the theorems. 

\subsection{Simulation results}
\begin{figure}[h]
\begin{tabular}{ll}
\includegraphics[scale=0.3]{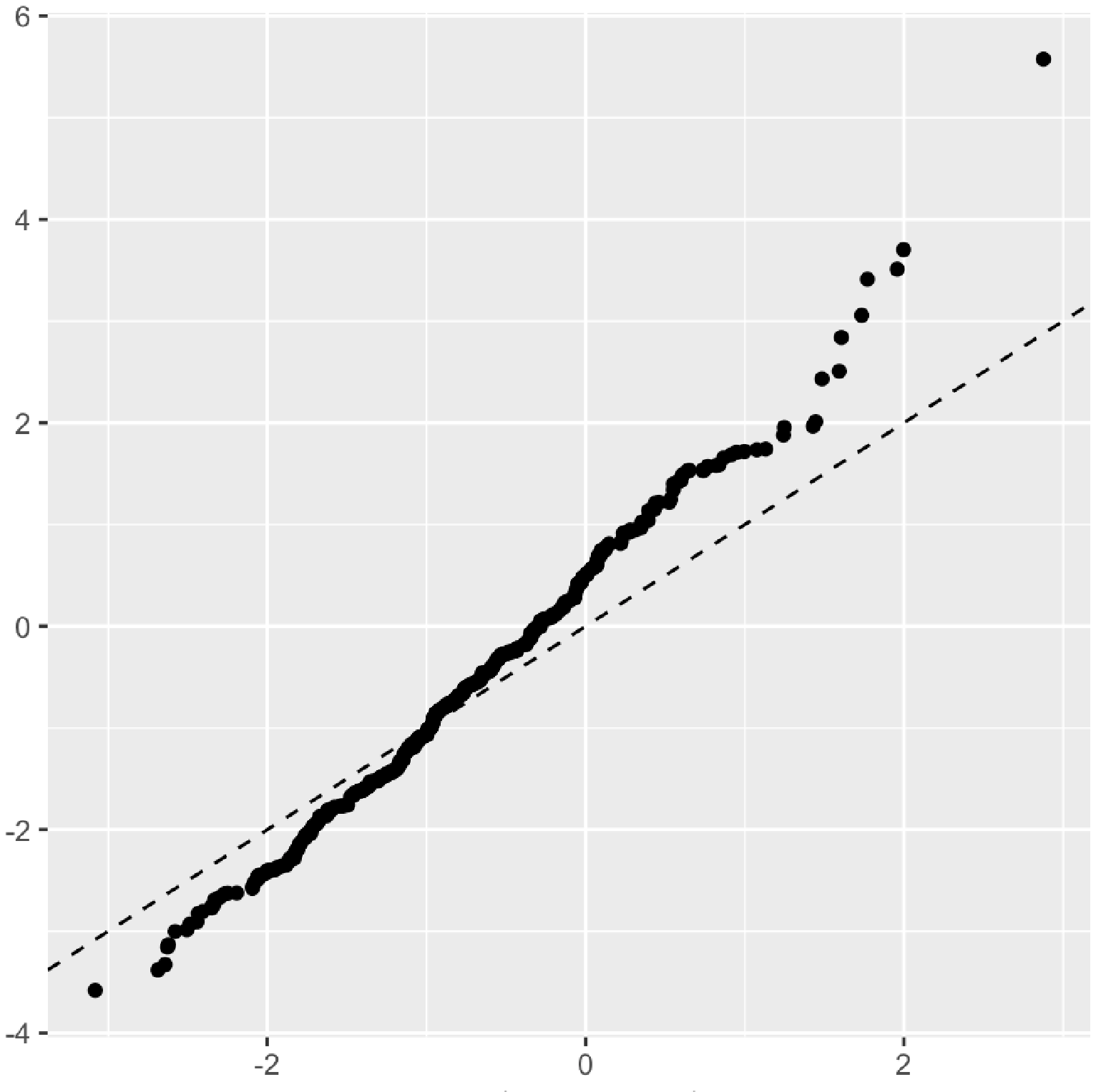}
&
\includegraphics[scale=0.3]{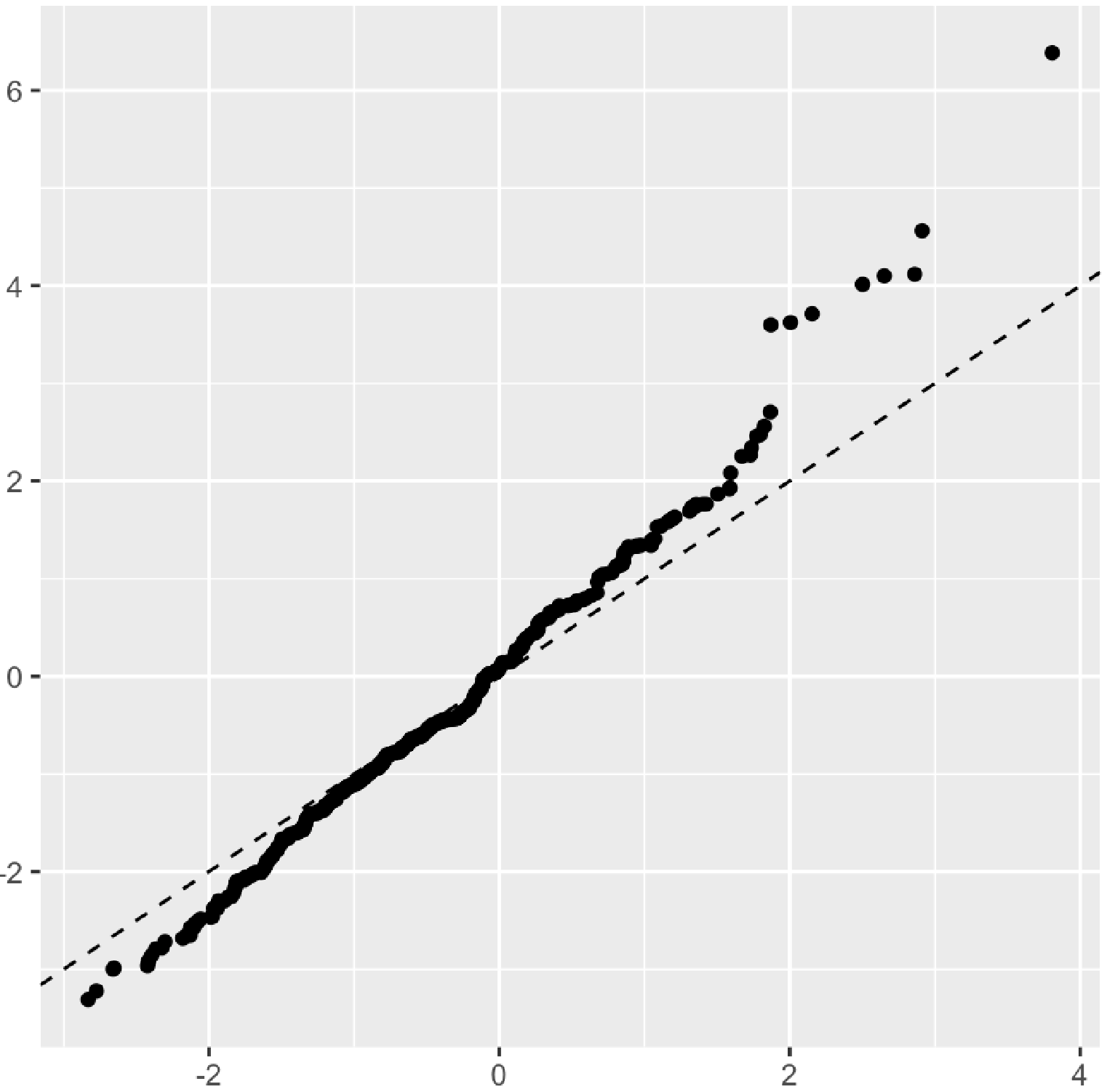}
\end{tabular}
\caption{$Q$-$Q$ plots corresponding to $(n, p)=(100,250)$ for the left picture and $(n, p)=(250, 500)$ for the right. The bold curves are empirical quantiles; the straight line is $y=x$.}
\label{Fig1}
\end{figure}

\begin{figure}[h]
\begin{tabular}{ll}
\includegraphics[width=4.7cm, height=4.6cm]{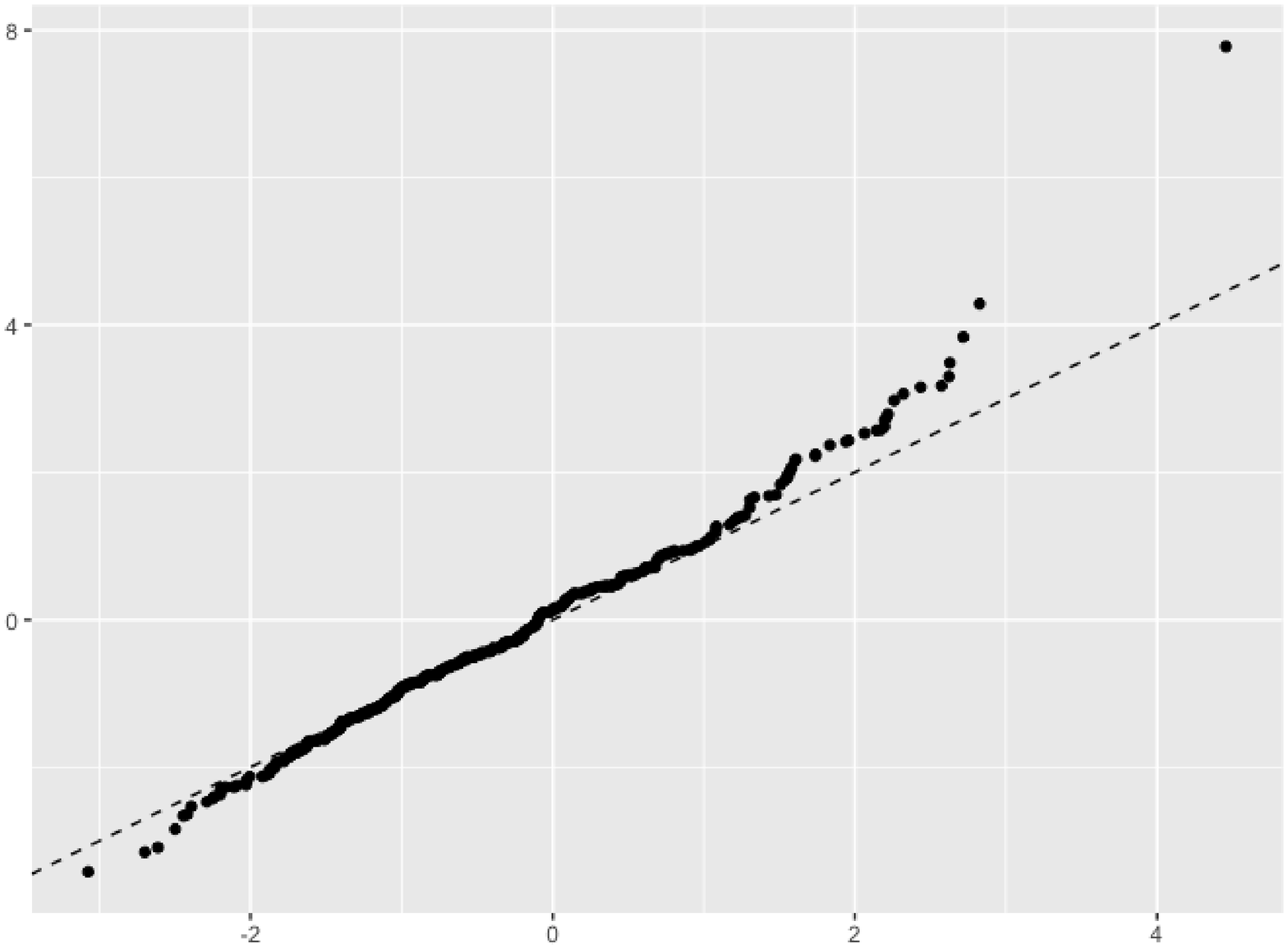}
&
\includegraphics[width=4.7cm, height=4.6cm]{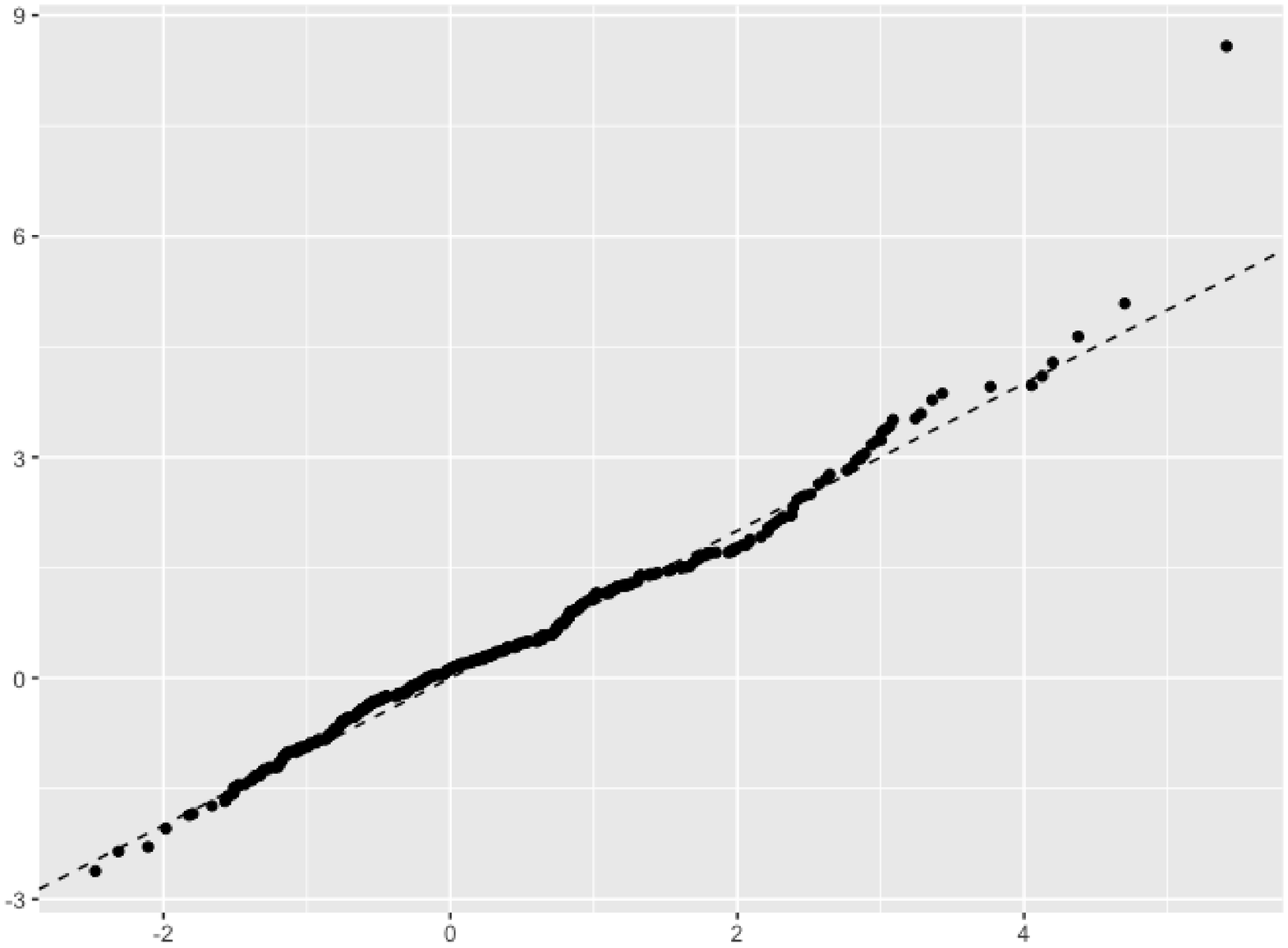}
\end{tabular}
\caption{$Q$-$Q$ plots corresponding to $(n, p)=(400,800)$ for the left picture and $(n, p)=(2000, 800)$ for the right. The bold curves are empirical quantiles; the straight line is $y=x$.}
\label{Fig2}
\end{figure}

\begin{figure}[h]
\begin{tabular}{ll}
\includegraphics[scale=0.38]{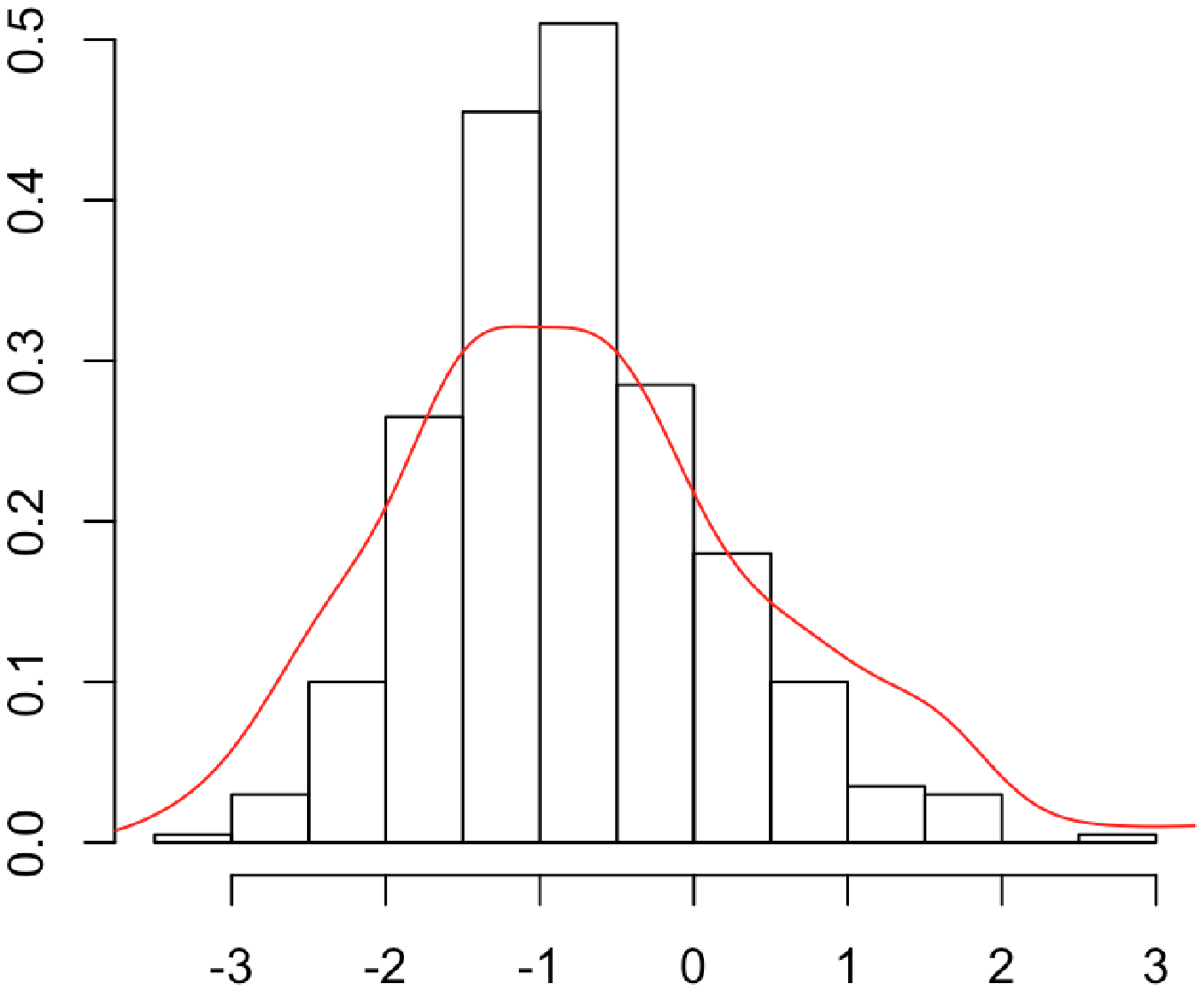}
&
\includegraphics[scale=0.38]{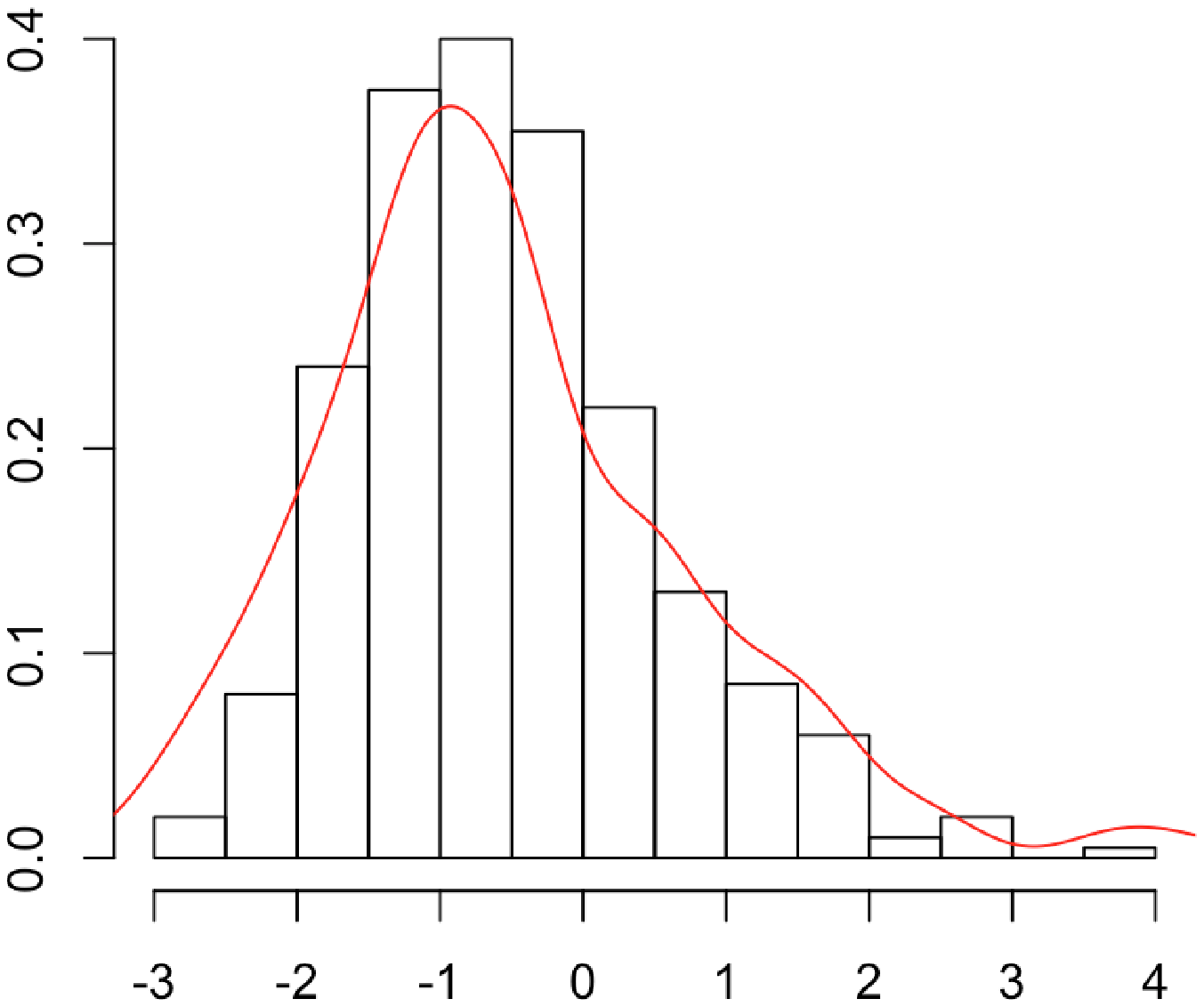}
\end{tabular}
\caption{Histograms of $W_n$ vs kernel density estimators of shifted  $W_n^*$ (red curves). Left: $W_n$ {\it vs} $W_n^{*}-1.2$ with  $(n, p)=(100, 250)$. Right: $W_n$ vs $W_n^{*}-1$ with  $(n, p)=(250, 500)$.}
\label{Fig3}
\end{figure}

\begin{figure}[h]
\begin{tabular}{ll}
\includegraphics[width=4.7cm, height=4.6cm]{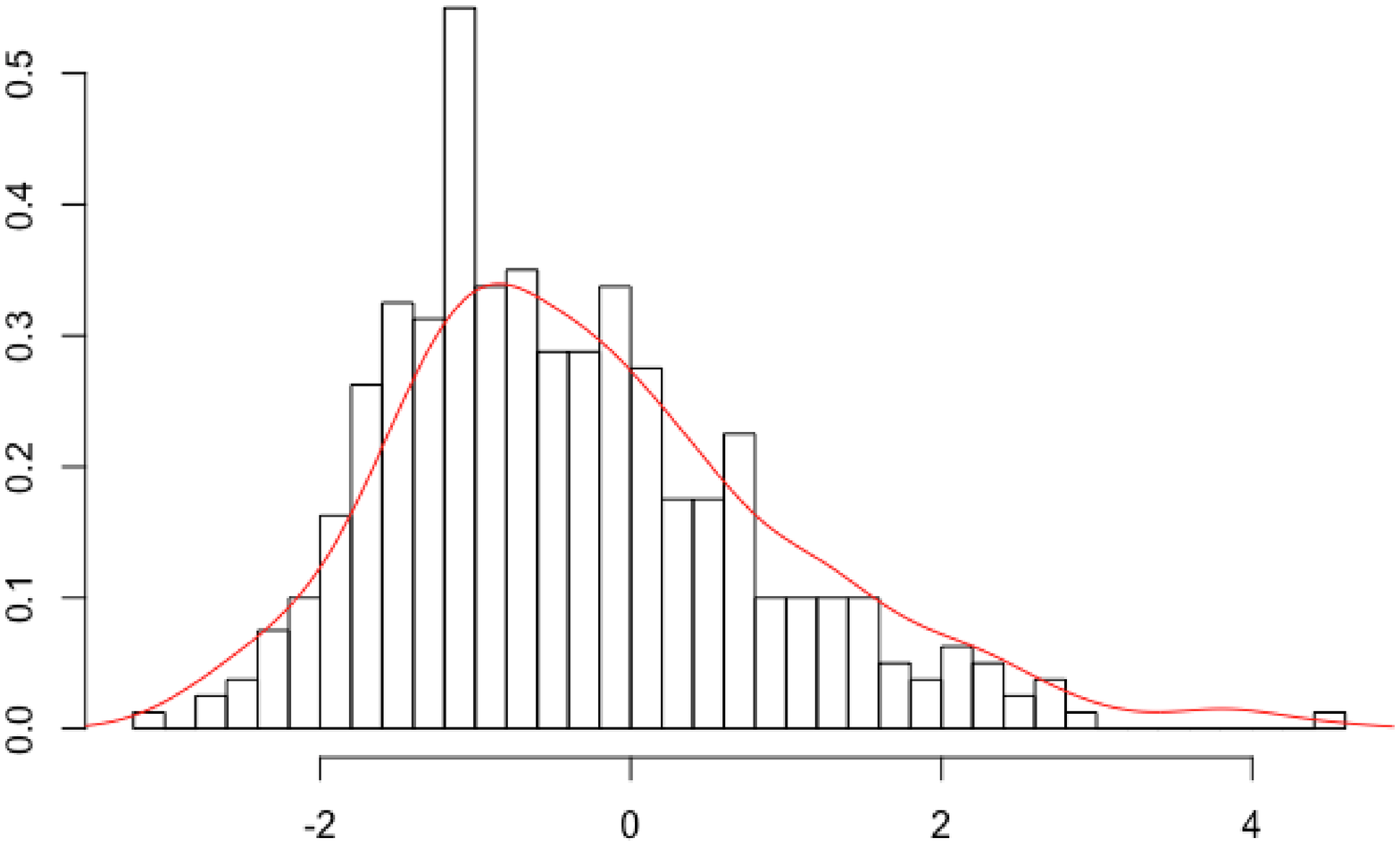}
&
\includegraphics[width=4.7cm, height=4.6cm]{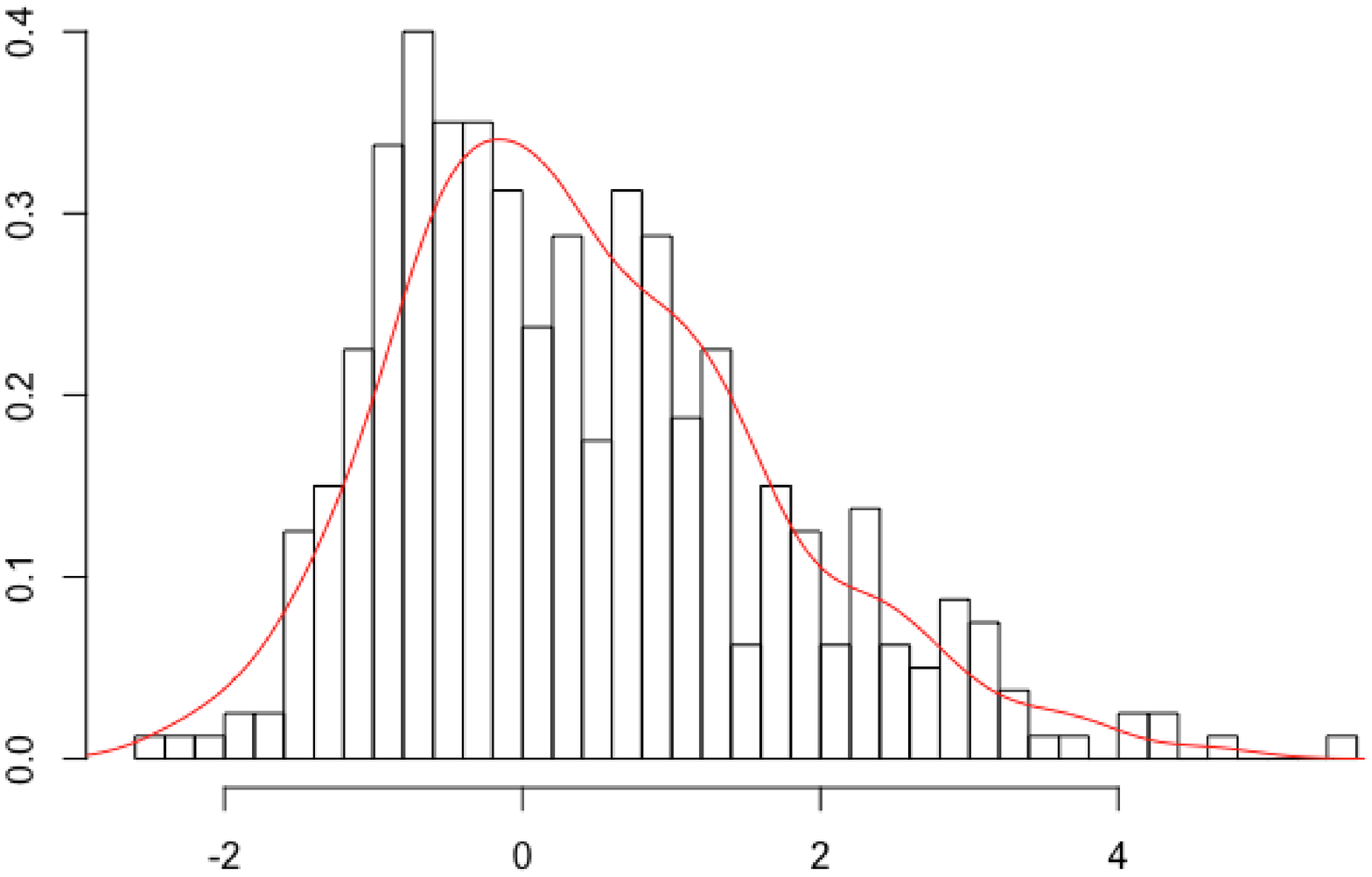}
\end{tabular}
\caption{Histograms of $W_n$ vs kernel density estimators of shifted  $W_n^*$ (red curves). Left:  $W_n$ {\it vs} $W_n^{*}-0.8$ with  $(n, p)=(400, 800)$. Right: $W_n$ vs $W_n^{*}$ with  $(n, p)=(2000, 800)$.}
\label{Fig4}
\end{figure}

We shall perform a Monte Carlo simulation to demonstrate the validity of our results. Let us describe the settings of our simulation: we have $m=400$ Monte Carlo iterations. Four sets of values for the pair $(n,p)$ are considered: $(n,p)=(100,250)$, $(n,p)=(200,500)$, $(n,p)=(400,800)$ and $(n,p)=(2000,800)$. We set $r=0.5$ for all the cases for simplicity. The Gaussian random variables $\bm{x}_1, \bm{x}_2,\cdots, \bm{x}_n$ have mean $\textbf{0}$ and $AR(1)$ covariance matrix with $r$ being the corresponding parameter. For each Monte Carlo iteration, we compute the value of 
\begin{align*}
W_n = \sqrt{2 \log p} \left(\frac{\sqrt{n} \mc{L}_n -r\sqrt{n}}{1-r^2}- \sqrt{2 \log p}+\frac{\log \log p+ \log 4\pi}{2\sqrt{2 \log p}} \right).
\end{align*}
Let $Z_1,Z_2,\cdots Z_p$ be i.i.d. standard normal random variables. Let $W_n^*$ be a random variables defined by
\begin{align*}
    W_n^{*}= \sqrt{2 \log p} \lb \max_{1 \leq i \leq p} Z_i - \sqrt{2 \log p} + \frac{\log \log p + \log(4 \pi)}{2 \sqrt{2 \log p}} \rb.
\end{align*}
It is well-known that $W_n^{*}$ converges weakly to a Gumbel distribution with the same cdf $F(x)=e^{-e^{-x}}$ appeared in Theorem \ref{slowdecay}. The convergence speed of $W_n^{*}$ is not fast. 
The proof of Theorem 1 suggests that $W_n$ and $W_n^*$ are close in distribution with difference being (roughly) of order $O_{\mb{P}}(n^{-1/2}(\log p)^{3/2})$; see Theorem S.1 in the supplement Jiang and Pham \cite{JiangPham21} for more details. Based on these, instead of comparing $W_n$ with $F(x)=e^{-e^{-x}}$ directly, we compare  $W_n$ and $W_n^{*}$ via $Q$-$Q$ plots  in Figures \ref{Fig1} and \ref{Fig2} and histograms in Figures \ref{Fig3} and \ref{Fig4}, respectively. A $Q$-$Q$ plot gives an insight on  how similar two distributions are, based on two data sets drawn from two distributions, respectively. When the two distributions are identical, the two $\alpha$ empirical quantiles are close. As a result, the curve formed by the pairs of $\alpha$ empirical quantiles for each $\alpha\in (0,1)$, which is called the $Q$-$Q$ plot, is close to the line $y=x$.

From Figures \ref{Fig1} and \ref{Fig2}, we see a bias between the bold curve and the line $y=x$ when $n$ and $p$ are small. This bias term is not surprising, as the proof of  Theorem \ref{slowdecay} suggests,  $\ds{W_n + O_{\mb{P}}(n^{-1/2}(\log p)^{3/2})}$ and $\ds{W_n^{*}}$ are close in distribution. The bias term $O_{\mb{P}}(n^{-1/2}(\log p)^{3/2})$ vanishes as both $n$ and $p$ are large. This can be seen from the second picture in Figure \ref{Fig2}.

Next, let us see the histograms in Figures \ref{Fig3} and \ref{Fig4}. Each histogram  describes the empirical distribution of $W_n$ versus a red curve, the kernel density estimators of the density of $W_n^*-c$. The normalizing constant ``$c$" in each figure reflects the bias term $O_{\mb{P}}(n^{-1/2}(\log p)^{3/2})$ mentioned above. The second picture from Figure \ref{Fig4} indicates that there is no obvious bias between the histogram and the red curve as $n$ and $p$ are larger. This is consistent with the $Q$-$Q$ plots considered earlier.

\section{Concluding remarks} \label{discussion}
In this paper, we investigate the largest Pearson  correlation coefficients of samples generated from ultra high-dimensional Gaussian populations. We study the problem under two scenarios: when the covariance matrix has an autoregressive covariance structure specified in (\ref{cov1}) and when it has a Toeplitz covariance structure as specified in (\ref{cov2}). Under the assumption that $\log p = o(n^{C})$ for some constant $C>0$ and $p=p_n \rightarrow \infty$, the asymptotic distributions of the largest Pearson correlation coefficients are obtained in each situation. In the case of $AR(1)$ model, the limiting distribution is Gumbel with different scaling   depending on how fast $r=r_n$ decays and the limiting distribution is the maximum of two independent Gumbel random variables at the critical case. In the case of the Toeplitz covariance structure, the limiting distribution is Gumbel if $r_n$ decays fast and it is a mixture of Gumbel and the standard normal distribution if $r_n$ is moderately large.  Furthermore, a central limit theorem is derived  for large $r_n$.
We now make some remarks.
    
    2. The assumptions $\log p = o(n^{1/5})$ and $\log p = o(n^{1/7})$ appeared in our main theorems may not be  optimal. Our techniques employed here require the two conditions. We conjecture that all the results hold as long as $\log p = o(n^{1/3})$. This improvement might be possible by carefully investigating the regularity/log-concavity of the distribution $S_k$ in the proof of Theorem \ref{slowdecay}. 
    
    3. In Theorem \ref{fastdecay}, it is natural to ask what would the limiting distribution be when $r$ converges to $1$ at a faster rate. It is possible that a central limit theorem holds. One would need a new technique to handle this case. Our techniques here for studying extreme values are not effective for random observations with strong dependence.  
    
    4.  We believe the classical technique involving the Poisson approximation can still be adapted to reproduce Theorems \ref{slowdecay}, \ref{large} and  \ref{fastdecay}. This is indeed our original approach to the problem. However, the technicalities are much more complicated and the argument is  quite lengthy. One advantage of this approach over the approach employed in this paper is the optimal dependence between $p$ and $n$. It is possible that the optimal condition $\log p = o(n^{1/3})$ can be achieved by using the method of Poisson approximation.
    
    5.  In Theorem \ref{Toeplitz}, an important question that needs to be addressed is the case $d \geq 2$ or more generally, when $\Sigma_n$ is allowed to change with $n$ and $d$ growing proportionally to $p$. When $d \geq 2$, the loss of stationarity  of the Gaussian field $Z_{i,i+s}^{n}$ (see the precise definition in proof of Theorem \ref{Toeplitz} below) prevents the uses of many classical techniques and also makes the covariance structure hard to understand. Another interesting direction is to study whether the condition $r_n \to 0$ in Theorem \ref{Toeplitz} can be dropped. We have seen in Example 4 that Theorem \ref{Toeplitz}  still holds even if $r_n$ does not converges to $0$.

\section{Proof of main results} \label{sec6}

Let us briefly outline the proofs. A more detailed explanation and heuristic derivations are given in the supplement Jiang and Pham \cite{JiangPham21}.

To prove Theorems \ref{slowdecay} and \ref{large}, the main tools are the high dimensional central limit theorems (Chernozukov et al.   \cite{Kato17, Kato20, Kato13, Kato19}), Bentkus \cite{Bentkus}, Gotze \cite{Gotze}, Portnoy \cite{Portnoy} and  Koike \cite{Koike}) and the Lipchitz concentration properties of Gaussian distributions. We first linearize the statistics $\hat{\rho}_{ij}$'s by using a Taylor expansion. Next, by carefully analyzing the error terms and separating out the main contribution terms, we show that the maximum $\mc{L}_n$ is close in distribution to the maximum of a Gaussian field with a non-constant mean function. We then combine the Lipchitz concentration properties of Gaussian distributions and classical techniques in extreme values theory to deduce the limiting distribution.

To prove Theorem \ref{fastdecay}, we first show that the maximum $\mc{L}_n$ is attained on the first sub-diagonal with probability going to one and then complete the proof by employing a similar mechanism used in the proofs of Theorems \ref{slowdecay} and \ref{large}.
Regarding Theorem \ref{Toeplitz}, the proof is similar to that of Theorem \ref{slowdecay} but with some appropriate adaptation. The main difference is that we have to deal with maximum of a non-stationary triangular array and the individual random variables are not weakly correlated. We overcome this difficulty by adapting the well-known normal comparison lemma for the  normal approximation. 

 During the proofs of main theorems, we will need some auxiliary results. The proofs of all lemmas in this section are presented in  the supplement by Jiang and Pham \cite{JiangPham21}

\subsection{Proof of Theorem \ref{slowdecay}} Before presenting the proofs, we need a few technical lemmas.  Recall the Orlicz norm $\|\cdot\|_{\psi_q}$ defined in \eqref{haoshi}. 
\begin{lemma} \label{vershynin}
Let $X,Y$ be sub-Gaussian random variables. Then,
$$ \| XY- \mb{E} (XY) \|_{\psi_1} \leq C \| XY \|_{\psi_1}\leq C\|X\|_{\psi_2} \|Y\|_{\psi_2}$$
where  $C>0$ is absolute constant.
\end{lemma}
For a vector $x=(x_1,x_2,\cdots,x_p)^T \in \mathbb{R}^p$, let $\| x \|_{\infty}=\max_{1 \leq i \leq p} |x_i|$. The next lemma gives an upper bound on the $\|.\|_{\infty}$-norm of sum of i.i.d. vectors.
\begin{lemma} \label{bernstein}
Let $\bm{X}_i=(X_{i1},X_{i2},\cdots,X_{ip})^T \in \mathbb{R}^p$, $i=1,2,\cdots,n$ be i.i.d.  random vectors with mean $\bm{0}$ and that $\sup_{i,j}\|X_{ij}\|_{\psi_1}<C$ for some constant  $C$ free of $n$ and $p$. Assume $p=p_n \rightarrow \infty$ and define
      $\bm{Z}_n=(\bm{X}_1+ \bm{X}_2+\cdots+ \bm{X}_n)/\sqrt{n}.$ 
Then $\|\bm{Z}_n\|_{\infty}=O_{\mb{P}}(\sqrt{\log p}).$  
\end{lemma}

Let  $\textbf{x}_1, \textbf{x}_2, \cdots, \textbf{x}_n$ be a random sample from a  $p$-dimensional normal population distribution with mean vector $\bm{0}$ and covariance matrix $\bm{\Sigma}$. Write $\textbf{x}_k=(x_{k1},x_{k2},\cdots,x_{kp})^T \in \mathbb{R}^p$ for each $k$. Then the data matrix is given by $(\textbf{x}_1, \textbf{x}_2, \cdots, \textbf{x}_n)^T=(x_{ij})_{n\times p}$ with all rows being i.i.d. $N(\bm{0}, \bm{\Sigma})$-distributed random vectors. With this setting we have the following result.

\begin{lemma} \label{Eij}
Suppose $\log p = o(n)$ and all diagonal entries of $\bm{\Sigma}$ are equal to $1$. Set
    \begin{align*}
        L_i &= -n + \sum_{k=1}^{n} x_{ki}^2 , ~~~~~
        \epsilon_{n,i} = \Big(1+\frac{L_{i}}{n}\Big)^{-1/2} - 1 + \frac{L_i}{2n}, \\
       E_{ij} &=\frac{L_iL_{j}}{4n^2}-\frac{L_i \epsilon_{n,j}}{2n}-\frac{L_{j} \epsilon_{n,i}}{2n}+ \epsilon_{n,i} \epsilon_{n,j}+ \epsilon_{n,i} + \epsilon_{n,j}.
    \end{align*}
    Then
        $$ \max_{1\leq i<j \leq p} \left\{|E_{ij}| \cdot \Big| \frac{1}{\sqrt{n}} \sum_{k=1}^{n} x_{ki}x_{kj} \Big|\right\} =O_{\mb{P}}\lr \frac{\log p}{\sqrt{n}} \rr.$$
\end{lemma}

\begin{lemma}\label{Corr} 
Let $x_{ij}$ be as stated above Lemma \ref{Eij} and $\bm{\Sigma}=(r^{|i-j|})_{p\times p}$ with $0\leq r \leq 1$. Define 
$S(i,j) = x_{1i}x_{1j} - \frac{1}{2} r^{|i-j|} (x_{1i}^2+x_{1j}^2).$
Then $\mbox{Var}(S(i,j)) = (1-r^{2|j-i|})^2$, 
\begin{eqnarray*}
      \mb{E}\left[S(i,i+1) \cdot S(j,j+1) \right]= \frac{1}{2}r^{2(j-i)} (1-r^2)^2~~\mbox{and}~~
    \sup_{(i,j) \neq (k,l)} | \mb{E}[S(i,j) \cdot S(k,l)]|
   \leq Cr,
\end{eqnarray*}
where the indices in the supremum also satisfy $1\leq i<j\leq p$,  $1\leq k<l\leq p$ and $C$ is an absolute constant.
\end{lemma}

We also need some extreme-value results of a weakly dependent random variables below. 
\begin{lemma} \label{berman}
Let $\{r_m>0; m\geq 1\}$ be numbers such that 
$1-r_m \geq C/\log m$ for all large $m$ and some constant $C$ free of $m$. Define

$$\bm{A}_m= \lr r_m^{2|i-j|} \lr 1 - \frac{I_{\la |i-j| \geq 1 \ra}}{2} \rr \rr_{1 \leq i, j \leq m} =\left(\begin{array}{cccc}
1 & \frac{r_m^2}{2} & \cdots & \frac{r_m^{2(m-1)}}{2} \\
\frac{r_m^2}{2} & 1 & \cdots & \frac{r_m^{2(m-2)}}{2} \\
\vdots & \vdots & & \vdots \\
\frac{r_m^{2(m-1)}}{2} & \frac{r_m^{2(m-2)}}{2} & \cdots & 1
\end{array}\right).$$
Let $\la X_{m,k} \ra_{1 \leq k \leq m}$ be a triangular array such that $(X_{m,1}, X_{m,2},\cdots, X_{m,m})^T \sim N(0, \bm{A}_m)$ for each $m\geq 1$.
Put $M_m=\max_{k \leq m} X_{m,k}$. Then, as $m\to\infty$, 
\begin{eqnarray*}
~~~~\mb{P} \lr \sqrt{2 \log m} \left( M_m -\sqrt{2\log m} + \frac{\log \log m + \log 4 \pi}{2 \sqrt{2 \log m}} \right) \leq x \rr \rightarrow e^{-e^{-x}}, ~~ x \in \mathbb{R}.
\end{eqnarray*}
\end{lemma}

\begin{lemma} \label{m2}
Let $\la  Z_{ij} \ra_{1 \leq i<j \leq p}$ be i.i.d. $N(0, \sigma_{ij}^2)$ with 
$\sigma_{ij}^2=(1-r^{2|i-j|})^2$ for some sequence $r=r_n\geq 0$ satisfying  $\limsup_{n\to\infty} r_n < 1$. Assume $\log p = o(n^{1/3})$ and 
$$\lim_{n \rar \infty} \frac{r \sqrt{n}}{\sqrt{\log p}} = L\in [0, \infty).$$
Denote 
$M^{(2)}_n = \max  \lb Z_{ij} + \sqrt{n} r^{|j-i|} \rb$,  where the maximum runs over all $1 \leq i<j \leq p$ with $j-i \geq 2$. Set $K_1=\frac{1}{2 \sqrt{2}}$.
Then,  as $n\to\infty$, 
$$\mb{P} \lr 2\sqrt{\log p} \lr M^{(2)}_n - 2\sqrt{\log p}+ \frac{\log \log p + \log 4 \pi}{4 \sqrt{\log p}} \rr \leq x \rr \rar e^{-K_1 e^{-x}},~~ x\in \mathbb{R}.$$
\end{lemma}

Now we are ready to prove Theorem \ref{slowdecay}. 

\noindent \textbf{Proof of Theorem \ref{slowdecay}}. We divide the proof into a few steps. 

 {\it Step 1: Linearization of the sample correlation coefficients.} We shall prove that 
 \begin{eqnarray}\label{linear}
 &&     \max_{1 \leq i < j \leq p} \sqrt{n} \hat{\rho}_{ij}\\  
 &=& \max_{1 \leq i < j \leq p} \la \sqrt{n} r^{|i-j|} + \frac{1}{\sqrt{n}} \sum_{k=1}^{n} \lb x_{ki}x_{kj} - \frac{r^{|i-j|}}{2}\left(x_{ki}^2+x_{kj}^2\right) \rb \ra + O_{\mb{P}}\Big(\frac{\log p}{\sqrt{n}}\Big) \nonumber
 \end{eqnarray}
 for any sequence $\{r_n\in [0, 1);\, n\geq 1\}$. 
 To see this, recall \eqref{wangshi} that
 \begin{eqnarray*}
\hat{\rho}_{i,j}=\frac{\sum_{k=1}^{n} x_{ki}x_{kj}}{\sqrt{\sum_{k=1}^{n} x_{ki}^2\sum_{k=1}^{n} x_{kj}^2}}.
\end{eqnarray*}
Set $L_i=\sum_{k=1}^{n} (x_{ki}^2-1)$ for $1\leq i \leq p$. Then  the denominator can be rewritten  as
 $$\frac{\sum_{k=1}^{n} x_{ki}^2}{n}=
 1+\frac{L_i}{n}.$$ 
 Inspired by the expansion $(1+x)^{-1/2}=1-x/2+O(x^2)$ as  $x\to 0$, recalling \eqref{wangshi}, we write
     \begin{align}
        \sqrt{n} \hat{\rho}_{ij}
        &=\left( \frac{1}{\sqrt{n}} \sum_{k=1}^{n} x_{ki}x_{kj}\right) \cdot\left( 1-\frac{L_i}{2n}+ \epsilon_{n,i} \right) \cdot\left( 1-\frac{L_{j}}{2n}+ \epsilon_{n,j} \right) \nonumber\\
        &= \left( \frac{1}{\sqrt{n}} \sum_{k=1}^{n} x_{ki}x_{kj}\right)\cdot \lb 2-\frac{1}{2n}\sum_{k=1}^{n} \left(x_{ki}^2+x_{kj}^2\right) + E_{ij} \rb \label{situation1}
    \end{align}
    where 
    $$
    \epsilon_{n,i}:=\Big(1+\frac{L_i}{n}\Big)^{-1/2}-1+\frac{L_i}{2n}
    $$
    and 
    $$E_{ij}:=\frac{L_iL_{j}}{4n^2}-\frac{L_i \epsilon_{n,j}}{2n}-\frac{L_{j} \epsilon_{n,i}}{2n}+ \epsilon_{n,i} \epsilon_{n,j}+ \epsilon_{n,i} + \epsilon_{n,j}.$$
   By Lemma \ref{Eij},
   \begin{eqnarray}\label{pingcheery} 
   \max_{i<j} \left\{|E_{ij}| \cdot \Big| \frac{1}{\sqrt{n}} \sum_{k=1}^{n} x_{ki}x_{kj} \Big|\right\} =O_{\mb{P}}\lr \frac{\log p}{\sqrt{n}} \rr. 
   \end{eqnarray}
  A further decomposition shows that
    \begin{align*}
       & \left( \frac{1}{\sqrt{n}} \sum_{k=1}^{n} x_{ki}x_{kj}\right) \cdot\lb 2-\frac{1}{2n}\sum_{k=1}^{n} \left(x_{ki}^2+x_{kj}^2\right)\rb\\
         = & \left[\frac{1}{\sqrt{n}} \sum_{k=1}^{n} \left(x_{ki}x_{kj} - r^{|i-j|}\right)\right]\cdot \lb 2-\frac{1}{2n}\sum_{k=1}^{n} \left(x_{ki}^2+x_{kj}^2\right)\rb
         +\\
         &~~~~~~~~~~~~~~~~~~~~~~~~~~~~~~~~~~~~~\sqrt{n}r^{|i-j|} \lb 2-\frac{1}{2n}\sum_{k=1}^{n} \left(x_{ki}^2+x_{kj}^2\right) \rb\\
         = &\lb \frac{1}{\sqrt{n}} \sum_{k=1}^{n} \left(x_{ki}x_{kj} - r^{|i-j|}\right)\rb\cdot \lb 1-\frac{1}{2n}\sum_{k=1}^{n} \left(x_{ki}^2+x_{kj}^2\right)\rb + Y_{ij}\\
         =& V_{ij}+Y_{ij},
    \end{align*}
    where $Y_{ij}$ and $V_{ij}$ are defined by
    \begin{align*}
            Y_{ij} &= \sqrt{n}r^{|i-j|} + \frac{1}{\sqrt{n}} \lb \sum_{k=1}^{n} x_{ki}x_{kj} - \frac{r^{|i-j|}}{2}\left(x_{ki}^2+x_{kj}^2\right) \rb;\\
            V_{ij} &= \lb \frac{1}{\sqrt{n}} \sum_{k=1}^{n} \left(x_{ki}x_{kj} -\sqrt{n} r^{|i-j|}\right)\rb \cdot\lb 1-\frac{1}{2n}\sum_{k=1}^{n} \left(x_{ki}^2+x_{kj}^2\right)\rb.
\end{align*}
Write
\begin{eqnarray*}
V_{ij} = -\lb \frac{1}{\sqrt{n}} \sum_{k=1}^{n} \big(x_{ki}x_{kj} -\mb{E}(x_{ki}x_{kj})\big)\rb \cdot\lb \frac{1}{n}\sum_{k=1}^{n} \frac{1}{2}\left(x_{ki}^2+x_{ki}^2-2\right)\rb.
\end{eqnarray*}
By Lemma \ref{vershynin}, $\sup_{1\leq i < j \leq p}\|x_{ki}x_{kj}-\mb{E}(x_{ki}x_{kj})\|_{\psi_1} \leq K$ and  $\sup_{1\leq i < j \leq p}\|(x_{ki}^2 + x_{kj}^2-2)/2\|_{\psi_1} \leq K$ for an absolute constant $K>0$. Set $\bm{X}_k= \la x_{ki}x_{kj}- \mb{E}(x_{ki}x_{kj}) \ra_{1 \leq i<j \leq p}\in \mathbb{R}^{p(p-1)/2}$, we obtain from Lemma \ref{bernstein} that 
    \begin{eqnarray*}
   \max_{1 \leq i < j \leq p}\Big|\frac{1}{\sqrt{n}}  \sum_{k=1}^{n} \big(x_{ki}x_{kj} - \mb{E}(x_{ki}x_{kj})\big)\Big|=\Big\|\frac{1}{\sqrt{n}}\big(\bm{X}_1+\cdots + \bm{X}_n\big)\Big\|_{\infty}= O_{\mb{P}}\Big(\sqrt{\log p}\Big).
    \end{eqnarray*}
   Similarly,   
     \begin{eqnarray*}
   \max_{1 \leq i < j \leq p}\frac{1}{n}\left|\sum_{k=1}^{n} \frac{1}{2}\big(x_{ki}^2+x_{ki}^2-2\big)\right|=O_{\mb{P}} \lr \frac{\sqrt{\log p}}{\sqrt{n}}\rr.
    \end{eqnarray*}
 Then
    \begin{eqnarray*}
   \max_{1 \leq i < j \leq p} |V_{ij}| 
  &\leq &  \max_{1 \leq i < j \leq p}\left|\frac{1}{\sqrt{n}} \sum_{k=1}^{n} \big(x_{ki}x_{kj} - \mb{E}(x_{ki}x_{kj})\big)\right|\cdot
   \max_{1 \leq i < j \leq p}\frac{1}{n}\left|\sum_{k=1}^{n} \frac{1}{2}\big(x_{ki}^2+x_{ki}^2-2\big)\right|\\
  &=& O_{\mb{P}}\Big(\sqrt{\log p}\Big) \cdot O_{\mb{P}} \lr \frac{\sqrt{\log p}}{\sqrt{n}} \rr= O_{\mb{P}}\Big(\frac{\log p}{\sqrt{n}}\Big).
    \end{eqnarray*}
This together with \eqref{situation1} and \eqref{pingcheery} leads to that $\max_{1 \leq i < j \leq p} \sqrt{n} \hat{\rho}_{ij}$ is identical to 
    \begin{eqnarray*}
    \max_{1 \leq i < j \leq p} \la \sqrt{n} r^{|i-j|} + \frac{1}{\sqrt{n}} \sum_{k=1}^{n} \lb x_{ki}x_{kj} - \frac{r^{|i-j|}}{2}(x_{ki}^2+x_{kj}^2) \rb \ra + O_{\mb{P}}\Big(\frac{\log p}{\sqrt{n}}\Big).
    \end{eqnarray*}
    Note Lemmas \ref{vershynin},  \ref{bernstein} and  \ref{Eij} hold based on the marginal information of $x_{ij} \sim N(0, 1)$ for each $i, j$ but not the covariances among them, then the above holds for any $r=r_n\in [0, 1).$ Thus (\ref{linear}) is proved.

{\it Step 2: Gaussian approximation to the second maximum in (\ref{linear})}. 
    For each $k=1,2,..,n$, consider the following collection of random variables
    \begin{eqnarray*}
    \bm{S}_k:=\la {x_{ki}x_{kj} - \frac{1}{2}r^{|i-j|}(x_{ki}^2+x_{kj}^2)} \ra_{1 \leq i<j \leq p} \in \mb{R}^{d}
    \end{eqnarray*}
    where $d:=p(p-1)/2$.  From now on, we shall also use the notation $S_k(i,j)$ to denote the $(i, j)$-coordinate of $\bm{S}_k$, that is,
    \begin{eqnarray*}
    S_k(i,j)={x_{ki}x_{kj} - \frac{1}{2}r^{|i-j|}(x_{ki}^2+x_{kj}^2)}.
     \end{eqnarray*}
  Since $ \mb{E}(x_{ki}x_{kj})=r^{|i-j|}$, we see $\mb{E}S_k(i,j)=0$. By assumption, $\textbf{x}_1, \textbf{x}_2, \cdots, \textbf{x}_n$ are i.i.d. random variables, thus $\bm{S}_1, \bm{S}_2,\cdots, \bm{S}_n$ are i.i.d. random vectors with mean $\bm{0}$. By Lemma \ref{Corr}, $\mbox{Var}(S_{1}(i,j))=(1-r^{2|j-i|})^2$ for any $1\leq i<j \leq p$. Moreover, by Lemma \ref{vershynin}, one has
    $$\sup_{1\leq i<j \leq p} \| S_{1}(i,j)\|_{\psi_1} < C$$
    for  some absolute constant $C$.  Let $\bm{N}=\la  N_{ij}; 1\leq i<j \leq p \ra$ be  $p(p-1)/2$ random variables which are jointly Gaussian with the same covariance structure as that of $\bm{S}_1$, that is,
    \begin{eqnarray}\label{zhuazhu}
    \la  N_{ij}; 1\leq i<j \leq p \ra 
    &\overset{d}{=} & \la  S_1(i, j); 1\leq i<j \leq p \ra\\
    &=& 
    \la x_{1i}x_{1j} - \frac{1}{2}r^{|i-j|}(x_{1i}^2+x_{1j}^2)); 1\leq i<j \leq p\ra.\nonumber
    \end{eqnarray}
    Thanks to high-dimensional central limit theorem (see Theorem S.2 in the supplement \cite{JiangPham21} and also Koike \cite{Koike}),  there exists an absolute constant  $C$ such that 
    \begin{align} \label{gauss}
                 \sup_{t\in \mb{R}^{d}} \bigg| \mb{P} \lr \frac{\sum_{k=1}^{n} \bm{S}_k}{\sqrt{n}} \leq t \rr  - \mb{P} \lr  \bm{N} \leq t \rr \bigg| \leq C\frac{(\log p)^{5/6}}{n^{1/6}}
    \end{align}
     where $t=(t_{ij})_{1\leq i<j \leq p}\in \mb{R}^{d}$.  
     Take $t_{ij}=s-\sqrt{n}r^{|j-i|}$ in (\ref{gauss}) to see 
        \begin{eqnarray} \label{gauss2}
             \sup_{s \in \mb{R}} \bigg| &&\mb{P} \lr \max_{1 \leq i<j \leq p} \left\{\frac{\sum_{k=1}^{n} S_k(i,j)}{\sqrt{n}} + \sqrt{n}r^{|j-i|}\right\} \leq s \rr  -\\
             &&\mb{P} \lr \max_{1 \leq i<j \leq p}\left\{ N_{ij} +\sqrt{n}r^{|j-i|}\right\} 
     \leq s \rr \bigg| 
     \leq  C\frac{(\log p)^{5/6}}{n^{1/6}}.\nonumber
    \end{eqnarray}
    By notation, $\mc{L}_n =\max_{1\leq i<j \leq p} \hat{\rho}_{i,j}$ and 
    \begin{eqnarray*}
    && \max_{1 \leq i<j \leq p} \left\{\frac{\sum_{k=1}^{n} S_k(i,j)}{\sqrt{n}} + \sqrt{n}r^{|j-i|}\right\}\\
    &= & \max_{1 \leq i < j \leq p} \la \sqrt{n} r^{|i-j|} + \frac{1}{\sqrt{n}} \sum_{k=1}^{n} \lb x_{ki}x_{kj} - \frac{r^{|i-j|}}{2}(x_{ki}^2+x_{kj}^2) \rb \ra\\
    &=& \sqrt{n}\mc{L}_n+ O_{\mb{P}}\Big(\frac{\log p}{\sqrt{n}}\Big)
    \end{eqnarray*}
    due to \eqref{linear}. It then follows from \eqref{gauss2} that
    \begin{eqnarray} \label{gausssha}
             && \sup_{s \in \mb{R}} \bigg| \mb{P} \lr \sqrt{n}\mc{L}_n+ O_{\mb{P}}\Big(\frac{\log p}{\sqrt{n}}\Big) \leq s \rr  -
             \mb{P} \lr \max_{1 \leq i<j \leq p}\left\{ N_{ij} +\sqrt{n}r^{|j-i|}\right\} 
     \leq s \rr \bigg| \\
     &\leq &   C\frac{(\log p)^{5/6}}{n^{1/6}},\nonumber
    \end{eqnarray}
    which tends to $0$ by the assumption $\log p = o(n^{1/5})$. Also, $O_{\mb{P}}((\log p)/\sqrt{n})$ appeared in the first probability is negligible after the normalization of $\sqrt{n\log p}\,\mc{L}_n$ as stated in (i) and (ii) of Theorem \ref{slowdecay}. Therefore, the proof of the theorem is reduced to analyzing the asymptotic distribution of the maximum of the Gaussian field $\bm{N}:=\{  N_{ij} +\sqrt{n} r^{|j-i|};\, 1 \leq i<j \leq p\}$. Keep in mind $N_{ij}\sim N(0, \sigma^2_{ij})$ with  $\sigma^2_{ij}=\mbox{Var}(S_{1}(i,j))=(1-r^{2|j-i|})^2$ and the covariance $\mbox{Cov}(N_{ij}, N_{kl})$ is estimated in Lemma \ref{Corr}.  
     
{\it Step 3: Analysis of maximum of Gaussian field $\bm{N}$.} 
We first show that the Gaussian field $\bm{N}$ can be approximated by a collection of centered, independent Gaussian random variables with similar means. 
        To this end, let $\bm{Z}=\la Z_{ij} \ra_{1 \leq i<j \leq p}$ be independent normal random variables with  $\mb{E}(Z_{ij})=\mb{E}(N_{ij})=0$ and  $\mbox{Var}(Z_{ij})= \mbox{Var}(N_{ij})=(1-r^{2|j-i|})^2$. By using an error bound for normal approximation (see Theorem S.3 in the supplement \cite{JiangPham21} and \cite{kato14-anti} for detailed discussion), one has
         $$ \sup_{t \in \mb{R}^{d}} \big| \mb{P} \lr \bm{N} \leq t \rr  - \mb{P} \lr  \bm{Z} \leq t \rr \big| \leq C \sup_{(i,j) \neq (k,l)} \bigg| \frac{\mb{E} \lb   N(i,j) \cdot N(k,l)  \rb}{(1-r^{2(j-i)})\cdot (1-r^{2(l-k)})} \bigg|^{1/3} \cdot (\log p)^{2/3}.$$
         By assumption, $\limsup_{n\to\infty} r_n < 1$. Note $ \sup_{(i,j) \neq (k,l)} |\mb{E}\lb N(i,j)\cdot N(k,l) \rb |=O(r)$ from Lemma \ref{Corr} and \eqref{zhuazhu}. Now choose $t=(t_{ij})_{1\leq i<j \leq p}$ with $t_{ij}=s-\sqrt{n}r^{|j-i|}$ to obtain
         \begin{eqnarray} \label{gauss3}
          &&   \sup_{s \in \mb{R}} \bigg| \mb{P} \lr \max_{1 \leq i<j \leq p} \left\{N_{ij} +\sqrt{n}r^{|j-i|} \right\}\leq s \rr  - \mb{P} \lr  \max_{1 \leq i<j \leq p} \left\{ Z_{ij} +\sqrt{n}r^{|j-i|}\right\} \leq s \rr \bigg| \\
          &\leq & O(r^{1/3}(\log p)^{2/3})=O\left(\frac{(\log p)^{5/6}}{n^{1/6}}\right) \to 0\nonumber
         \end{eqnarray}
         by the assumption $\log p = o(n^{1/5})$ and therefore, the problem is reduced  to studying 
         \begin{eqnarray}\label{weishaai}
         M_n: = \max_{1 \leq i<j \leq p} \big\{Z_{ij} +\sqrt{n}r^{|j-i|}\big\}.
         \end{eqnarray}
         The advantage of this reduction is that $\{Z_{ij};\, 1\leq i<j \leq p\}$ are independent. Recall $L=\lim_{n \rightarrow \infty} r\sqrt{n}(\log p)^{-1/2}\in [0, \infty].$ To complete the proof, we will consider three cases: $ 0 \leq L < 2 - \sqrt{2}$, $2-\sqrt{2} < L < \infty$ and $L=\infty$ separately in the following. 
         
          {\it Case 1: $ 0 \leq L < 2 - \sqrt{2}$}.
         Define  
         \begin{align}
             M^{(1)}_n = \max_{j-i=1} \la Z_{ij} + r \sqrt{n} \ra \ \ \ \mbox{and}\ \ \ 
             M^{(2)}_n = \max_{j-i \geq 2}  \la Z_{ij} + \sqrt{n} r^{|j-i|} \ra. \label{mdef1}
         \end{align}
        Obviously, $M_n=\max\{M^{(1)}_n, M^{(2)}_n\}$. 
        We shall prove that $M^{(2)}_n > M^{(1)}_n$ with probability tending to $1$. To achieve this, we first estimate $\mb{E} M^{(1)}_n$ and $\mb{E} M^{(2)}_n$. Set $Z^{*}_{ij}=Z_{ij}/(1-r^2)$. Then $\{Z^{*}_{ij};\, 1\leq i<j \leq p\}$ are i.i.d. standard normals. One can rewrite $\mb{E} M^{(1)}_n$ as
        \begin{eqnarray}
        \mb{E} M^{(1)}_n = r \sqrt{n} + (1-r^2)\cdot \mb{E} \max_{j-i=1} Z^{*}_{ij}
        &=& r\sqrt{n} + \sqrt{2 \log (p-1)}(1+o(1))\nonumber\\
        &=& (L+\sqrt{2})\sqrt{\log p}\,(1+o(1)) \label{daren}
        \end{eqnarray}
        where we use a formula on the expected value of the maximum of $p$ i.i.d. standard normals; see, for example,  Exercise 2.11 from \cite{MW2019}. Note that
        $$\max_{j-i \geq 2} \left\{r^{|i-j|} \sqrt{n}\right\} \leq r^2\sqrt{n}\leq \frac{\log p}{\sqrt{n}} \rightarrow 0.$$
        Moreover, form the expression that $M^{(2)}_n = \max_{j-i \geq 2}  \{ (1-r^2)Z^*_{ij} + \sqrt{n} r^{|j-i|}\}$ we see that
        \begin{eqnarray*}
        \Big|M^{(2)}_n -\max_{j-i \geq 2} Z^*_{ij}\Big|
        \leq r^2\cdot \max_{j-i \geq 2}  |Z^*_{ij}| + \sqrt{n} r^2.
        \end{eqnarray*}
        Observe both maxima are taken over $\frac{1}{2}p^2(1+o(1))$ i.i.d. standard normals. Thus $E\max_{j-i \geq 2} Z^*_{ij}=\sqrt{2\log ((p^2/2)(1+o(1)))}=2\sqrt{\log p} (1+o(1))$ and 
        \begin{eqnarray*}
        \Big| \mb{E}M^{(2)}_n - \mb{E}\max_{j-i \geq 2} Z^*_{ij}\Big|
        & \leq & r^2\cdot \mb{E}\max_{j-i \geq 2}  |Z^*_{ij}| +
        \sqrt{n} r^2\\
        & \leq &  O\Big(\frac{\log p}{n}\Big) \cdot 3\sqrt{\log p} + O\left(\frac{\log p}{\sqrt{n}}\right)=O\left(\frac{\log p}{\sqrt{n}}\right)\to 0
        \end{eqnarray*}
        by the assumption $\log p = o(n^{1/5})$. 
        Consequently,
        \begin{eqnarray}\label{SQSX}
        \mb{E} M^{(2)}_n = 2 \sqrt{\log p}(1+o(1)).
        \end{eqnarray}
        Comparing this with \eqref{daren} and using the assumption $0 \leq L < 2 - \sqrt{2}$, we know there exists  $\delta>0$ small enough such that 
        $$
        \mb{E} M^{(1)}_n + \delta \sqrt{\log p}<\mb{E} M^{(2)}_n - \delta \sqrt{\log p}
        $$
        for large $n$. 
        By the Lipschitz concentration inequality for Gaussian distributions (see \cite{MW2019}, Example 2.29), one has 
        \begin{eqnarray*}
           &&\mb{P}(M^{(1)}_n - \mb{E} M^{(1)}_n > \delta \sqrt{\log p} ) \leq \exp\Big\{-\frac{C\delta^2 \log p}{2}\Big\}; \\
            &&\mb{P}(M^{(2)}_n - \mb{E} M^{(2)}_n \leq -\delta \sqrt{\log p} )  \leq \exp\Big\{-\frac{C\delta^2 \log p}{2}\Big\}.
        \end{eqnarray*}
        Therefore, with probability at least $1-2\exp\{-\frac{C\delta^2 \log p}{2}\}$ and for sufficiently large $n$, one has $M^{(1)}_n < M^{(2)}_n$.
        Hence, with probability at least $1-2e^{-\frac{C\delta^2 \log p}{2}}$, it holds that $M_n = M^{(2)}_n$.
        This fact together with  Lemma \ref{m2} yields
        \begin{eqnarray}\label{rezhong}
        \mb{P} \lr 2\sqrt{\log p} \lr M_n - 2\sqrt{\log p}+ \frac{\log \log p + \log 4 \pi}{4 \sqrt{\log p}} \rr \leq x \rr \rar e^{-K_1 e^{-x}},
        \end{eqnarray}
         where $K_1=\frac{1}{2 \sqrt{2}}$. Combining \eqref{gausssha},  \eqref{gauss3} with \eqref{weishaai}
         the above convergence, we arrive at
         $$\mb{P} \lr 2\sqrt{\log p} \lr \sqrt{n} \mc{L}_n - 2\sqrt{\log p}+ \frac{\log \log p + \log 4 \pi}{4 \sqrt{\log p}}\rr +O_{\mb{P}}\lr\sqrt{\log p}\cdot \frac{\log p}{\sqrt{n}}\rr \leq x \rr $$ 
         converges to $e^{-K_1 e^{-x}}$ for every $x\in \mathbb{R}$, 
where $K_1=\frac{1}{2 \sqrt{2}}$. By the assumption $\log p =o(n^{1/5})$, the ``$O_{\mb{P}}$'' term in the above probability goes to $0$, then the Slutsky lemma implies  
 $$ 2\sqrt{\log p} \lr \sqrt{n} \mc{L}_n - 2\sqrt{\log p}+ \frac{\log \log p + \log 4 \pi}{4 \sqrt{\log p}}\rr $$ 
converges weakly to the distribution with cdf $e^{-K_1 e^{-x}}$. This concludes (i) of Theorem \ref{slowdecay}. 

{\it Case 2:  $2-\sqrt{2} < L < \infty$}. One can derive the limiting distribution by using a similar argument to the {\it Case 1} with a minor change. Recalling \eqref{mdef1}, we will show $M^{(1)}_n>M^{(2)}_n$ with probability converging to $1$, and hence $M_n = M^{(1)}_n$ with probability going to $1$.

         Recall \eqref{daren} and \eqref{SQSX}. Choose $\delta>0$ small enough such that 
         $$
        \mb{E} M^{(1)}_n - \delta \sqrt{\log p}> \mb{E} M^{(2)}_n +\delta \sqrt{\log p}
        $$
        for large $n$. 
         By the Lipchitz concentration inequality again,
        \begin{eqnarray*}
        &&    \mb{P}(M^{(1)}_n - \mb{E} M^{(1)}_n \leq -\delta \sqrt{\log p} ) \leq \exp\left\{-\frac{C\delta^2 \log p}{2}\right\}, \\
        &&   \mb{P}(M^{(2)}_n - \mb{E} M^{(2)}_n \geq \delta \sqrt{\log p} )  \leq \exp\left\{-\frac{C\delta^2 \log p}{2}\right\}.
        \end{eqnarray*}
        Therefore, with probability at least $1-2e^{-\frac{C\delta^2 \log p}{2}}$, one has $M^{(1)}_n > M^{(2)}_n$ as $n$ is sufficiently large, and hence 
        $M_n = \max\{M^{(1)}_n, M^{(12}_n\}=M^{(1)}_n$. Easily,  $M^{(1)}_n-r\sqrt{n}$ is the maximum of $(p-1)$ i.i.d. random variables of distribution $N(0, (1-r^2)^2)$. Hence, by a standard result on the maximum of standard normals (see, e.g., Theorem 1.5.3 from \cite{Leadbetter}) that
\begin{eqnarray}\label{qianjin}
\sqrt{2\log (p-1)} \lr \frac{M_n -r\sqrt{n}}{1-r^2} - \sqrt{2\log (p-1)}+ \frac{\log \log (p-1) + \log 4 \pi}{2 \sqrt{2\log (p-1)}} \rr
\end{eqnarray}
converges weakly to the distribution with cdf $e^{-e^{-x}}$. Since $\sqrt{\log (p-1)}-\sqrt{\log p}=O(1/(p\log p))$ and a similar estimate holds for $\log \log p$, by Slusky's lemma,
        \begin{eqnarray}\label{jiade}
        \sqrt{2\log p} \lr \frac{M_n -r\sqrt{n}}{1-r^2} - \sqrt{2\log p}+ \frac{\log \log p + \log 4 \pi}{2 \sqrt{2\log p}} \rr
        \end{eqnarray}
        converges weakly to $e^{-e^{-x}}$.
        By the same argument after \eqref{rezhong}, we arrive at
      $$\mb{P} \lr \sqrt{2\log p} \lr \frac{ \sqrt{n} \mc{L}_n  -r\sqrt{n}}{1-r^2} - \sqrt{2\log p}+ \frac{\log \log p + \log 4 \pi}{2 \sqrt{2\log p}} \rr \leq x \rr \rar e^{- e^{-x}}.$$
      Thus,  (ii) of Theorem \ref{slowdecay} for the case $L<\infty$ has been proved.

{\it Case 3:  $L = \infty$}. Under this assumption, $\lim_{n \rightarrow \infty} r\sqrt{n}(\log p)^{-1/2}=\infty.$ Recall {\it Step 2} and the identity in \eqref{zhuazhu} in particular, we know that $\{  N_{ij} +\sqrt{n} r^{|j-i|};\, 1 \leq i<j \leq p\}$ are jointly Gaussian random variables with $N_{ij}\sim N(0, \sigma^2_{ij})$,    $\sigma^2_{ij}=(1-r^{2|j-i|})^2$ 
 and $\mbox{Cov}(N_{ij}, N_{kl})$ is given in Lemma \ref{Corr}.  Define
         \begin{align*}
             Q^{1}_n = \max_{j-i=1} \left\{N_{ij} +r\sqrt{n}\right\}\ \ \ \ \ \mbox{and}\ \  \ \ \ \ 
             Q^{2}_n = \max_{j-i \geq 2} \left\{ N_{ij} + \sqrt{n}r^{|i-j|} \right\}.
         \end{align*}
         Then
         \begin{eqnarray}\label{sinian}
         \max_{1 \leq i<j \leq p}\left\{ N_{ij} +\sqrt{n}r^{|j-i|}\right\}=\max\left\{ Q^{1}_n,  Q^{2}_n  \right\}.
         \end{eqnarray}
         We shall show that, with probability going to $1$, 
         $$Q^1_n \geq \frac{(r+r^2)\sqrt{n}}{2} \geq Q^2_n.$$

        In fact, observe $(r^2-r)\sqrt{n}[2(1-r^2)]^{-1}\to -\infty$ since $\lim_{n \rightarrow \infty} r\sqrt{n}(\log p)^{-1/2}=\infty.$ Then 
         \begin{align}
                     \mb{P} \lr Q^{1}_n  \geq \frac{(r+r^2)\sqrt{n}}{2} \rr & \geq  \mb{P} \lr N_{12} + r\sqrt{n} \geq  \frac{(r+r^2)\sqrt{n}}{2} \rr \nonumber\\
                     & =\mb{P} \lr N(0, 1)  \geq  \frac{(r^2-r)\sqrt{n}}{2(1-r^2)} \rr \to 1. \label{kudai}
         \end{align}
         We now consider the event $\la Q^2_n \leq \frac{(r+r^2)\sqrt{n}}{2} \ra$. 
        By the union bound, one has that
         \begin{align}
              \mb{P} \lr Q^{2}_n  > \frac{(r+r^2)\sqrt{n}}{2} \rr 
              &=\mb{P} \lr \max_{j-i \geq 2} \left\{ N_{ij} + \sqrt{n}r^{|i-j|} \right\}  > \frac{(r+r^2)\sqrt{n}}{2} \rr \nonumber\\
              &\leq \sum_{j-i \geq 2} \mb{P} \lr N_{ij} \geq \frac{\sqrt{n}(r+r^2-2r^{|i-j|})}{2} \rr.\label{hongqi}
         \end{align}
         Evidently,
         $$
         \frac{\sqrt{n}(r+r^2-2r^{|i-j|})}{1-r^{2|j-i|}} \geq \sqrt{n}(r-r^2)
         $$
         for $j-i \geq 2$. Since $N_{ij}/(1-r^{2|j-i|}) \sim N(0, 1)$, the last sum in \eqref{hongqi} is bounded by
         \begin{eqnarray}
          p^2\cdot
         \mb{P} \lr N(0, 1)\geq \frac{1}{2}\sqrt{n}(r-r^2)\rr \leq  p^2\cdot \exp\left\{-\frac{1}{8}nr^2(1-r)^2\right\}\label{tonghua}
         \end{eqnarray}
         thanks to the well-known Gaussian tail bound $P(N(0, 1)\geq x) \leq e^{-x^2/2}$ for all $x\geq 0$. By  assumption, we know $\lim_{n \rightarrow \infty} r\sqrt{n}(\log p)^{-1/2}=\infty$ and $\limsup_{n\to\infty} r_n < 1$. These imply that $nr^2(1-r)^2(\log p)^{-1}\to \infty$. Consequently, the last term in \eqref{tonghua} goes to $0$. It follows that $\mb{P} \lr Q^{2}_n  \leq (r+r^2)\sqrt{n}/2 \rr \to 1$. This and \eqref{kudai} yield $\mb{P} \lr Q^{1}_n>Q^{2}_n\rr \to 1$.  By \eqref{sinian},
         $$
         \max_{1\leq i<j \leq p}\left\{N_{ij} +\sqrt{n}r^{|j-i|}\right\}=Q_n^1=\max_{1\leq i \leq p-1} \left\{N_{i,i+1} \right\}+r\sqrt{n}
         $$
         with probability going to $1$. Here and below, to not cause any confusion, with a bit abuse of notation, we use ``$N_{i,j}$'' for ``$N_{ij}$'' if necessary. Recall $N_{i,i+1}\sim N(0, (1-r^2)^2)$. Then
         \begin{eqnarray}\label{Beasly}
         \max_{1\leq i<j \leq p}\left\{N_{ij} +\sqrt{n}r^{|j-i|}\right\}
         =(1-r^2)\cdot \max_{1\leq i \leq p-1} U_i+r\sqrt{n}
         \end{eqnarray}
        with probability going to $1$ and $U_i:=(1-r^2)^{-1}N_{i,i+1}$. Notice each of  $\{U_i;\, 1\leq i \leq p-1\}$ has the standard normal distribution and 
        $$
        \mbox{Cov}(U_i, U_j)=\frac{1}{2}r^{2(j-i)},~~~ 1\leq i<j\leq p-1
        $$
        by Lemma  \ref{Corr} (match \eqref{zhuazhu} and notation $S(i, j)$ in the lemma).  Then $(U_1, \cdots, U_{p-1})^T\sim N(0, \bm{A}_{p-1})$ where $\bm{A}_{p-1}$ is given in Lemma \ref{berman}. By this lemma and \eqref{Beasly},
        \begin{eqnarray*}
        &&\sqrt{2\log (p-1)}\left(\frac{1}{1-r^2}\max_{1\leq i<j \leq p}\left\{N_{ij} +\sqrt{n}r^{|j-i|}-r\sqrt{n}\right\} \right) \nonumber\\
        &- &2\log (p-1)+\frac{1}{2}\left[\log\log (p-1) +\log (4\pi)\right] 
        \end{eqnarray*}
         converges weakly to a distribution with cdf $e^{-e^{-x}}$ for $x \in \mathbb{R}.$ Multiplying the above by $\sqrt{2\log p}/\sqrt{2\log (p-1)}$, noticing $\log (p-1)=\log p+o(1)$ and $\log \log (p-1)=\log \log p +o(1)$, we obtain from the Slutsky lemma that
          \begin{eqnarray*}
          \sqrt{2\log p}\left(\frac{1}{1-r^2}\max_{1\leq i<j \leq p}\left\{N_{ij} +\sqrt{n}r^{|j-i|}-r\sqrt{n}\right\} \right)
        - 2\log p+\frac{1}{2}\left[\log\log p +\log (4\pi)\right]
        \end{eqnarray*}
         converges weakly to the distribution with cdf $e^{-e^{-x}}$. From \eqref{gausssha} we see that
         \begin{eqnarray*}
          \sqrt{2\log p}\cdot\frac{\sqrt{n}\mathcal{L}_n-r\sqrt{n}+O(n^{-1/2}\log p)}{1-r^2}-2\log p+\frac{1}{2}\left[\log\log p +\log (4\pi)\right]
        \end{eqnarray*}
         converges weakly to  $e^{-e^{-x}}$. We get  (ii) of Theorem \ref{slowdecay} for $L=\infty$ by Slutsky's lemma. \hfill$\square$

\subsection{Proof of Theorem \ref{large}}
The idea of the proof is similar to that of Theorem \ref{slowdecay}. We will make some modifications and carry out a finer analysis based on the convergence speed of  $r \sqrt{n/\log p}\to 2-\sqrt{2}$. We will continue to use the notation in the proof of Theorem \ref{slowdecay}.

Review {\it Steps 1, 2} and the case $0<L<2-\sqrt{2}$ in  {\it Step 3} in the proof of Theorem \ref{slowdecay}, the argument is still valid since $L$ is finite. In particular, we have from \eqref{gausssha} and \eqref{gauss3} that
\begin{eqnarray} \label{earth2}
             \sup_{s \in \mb{R}} \bigg| \mb{P} \lr \sqrt{n}\mc{L}_n+ O_{\mb{P}}\Big(\frac{\log p}{\sqrt{n}}\Big) \leq s \rr  - \mb{P} \lr  M_n \leq s \rr \bigg| 
          = O\left(\frac{(\log p)^{5/6}}{n^{1/6}}\right) \to 0,\nonumber
    \end{eqnarray}
where 
$
M_n:=\max_{1 \leq i<j \leq p} \left\{ Z_{ij} +\sqrt{n}r^{|j-i|}\right\}
$
and $\la Z_{ij} \ra_{1 \leq i<j \leq p}$ are independent normal random variables with  $\mb{E}(Z_{ij})=0$ and  $Var(Z_{ij})=(1-r^{2|j-i|})^2$. Recall the notation
\begin{align*}
             M^{(1)}_n = \max_{j-i=1} \la Z_{ij} + r \sqrt{n} \ra, \ \ \ \ \ \ 
             M^{(2)}_n = \max_{j-i \geq 2}  \la Z_{ij} + \sqrt{n} r^{|j-i|} \ra,
         \end{align*}
         where the indices $i,j$ in the above maxima also satisfy $1\leq i<j \leq p$. 
         Evidently,  $M^{(1)}_n$ and $M^{(2)}_n$ are independent and 
         $
        M_n:=\max\{M^{(1)}_n, M^{(2)}_n\}.
        $
The purpose is to derive the asymptotic distribution of $\mc{L}_n$. Through \eqref{earth2} the problem is reduced to the study of  the limits of $M^{(1)}_n$ and $M^{(2)}_n$, respectively. We investigate these next.

Similar to the discussions in \eqref{qianjin} and \eqref{jiade}, we have 
\begin{eqnarray*}
        J_n:=\sqrt{2\log p} \lr \frac{M^{(1)}_n -r\sqrt{n}}{1-r^2} - \sqrt{2\log p}+ \frac{\log \log p + \log 4 \pi}{2 \sqrt{2\log p}} \rr \to G
        \end{eqnarray*}
        weakly,  where $G$ is the Gumbel distribution with cdf  $e^{-e^{-x}}$. In particular,
        \begin{eqnarray}\label{tata}
        \lim_{n\to\infty} \mb{P}(J_n\leq y_n)=
        \begin{cases}
        1, ~~\text{if $y_n\to\infty$};\\
        \exp(-e^{-y}), ~~\text{if $y_n\to y$};\\
        0, ~~\text{if $y_n\to-\infty$}.
        \end{cases}
        \end{eqnarray}
        Fix $x\in \mb{R}$, define $t=t(x)$ by
\begin{eqnarray}\label{toujiang}
t= 2\sqrt{\log p} +\frac{x}{2 \sqrt{\log p}}  - \frac{\log \log p + \log (4 \pi)}{4 \sqrt{\log p}}.
\end{eqnarray}
From {\it Step 3} of the proof of Theorem \ref{slowdecay}, one has
\begin{eqnarray}\label{pianzi}
\mb{P}(M^{(2)}_n \leq t) \rar e^{-K_1 e^{-x}}
\end{eqnarray}
where $K_1 = (2\sqrt{2})^{-1}$. To analyze the asymptotic property of $\mb{P}(M^{(1)}_n \leq t)$, we rewrite 
\begin{eqnarray}\label{mangong}
    \mb{P}(M^{(1)}_n \leq t) 
    &=& \mb{P} \lr \sqrt{2\log p} \lr \frac{M^{(1)}_n -r\sqrt{n}}{1-r^2} - \sqrt{2\log p}+ \frac{\log \log p + \log (4 \pi)}{2 \sqrt{2\log p}} \rr \leq t_1 \rr\nonumber\\
    &=& \mb{P} \lr J_n \leq t_1\rr
\end{eqnarray}
where 
\begin{align}\label{whysha}
    t_1:&= \sqrt{2\log p} \lr \frac{t -r\sqrt{n}}{1-r^2} - \sqrt{2\log p}+ \frac{\log \log p + \log (4 \pi)}{2 \sqrt{2\log p}} \rr \nonumber\\
    & = \sqrt{2\log p} \cdot \frac{t -r\sqrt{n}}{1-r^2} - 2 \log p +\frac{\log \log p + \log 4\pi}{2} \nonumber\\
    &= \sqrt{2 \log p}\,(t-r\sqrt{n}) - 2 \log p +\frac{\log \log p + \log 4\pi}{2} +o(1).
\end{align}
In the above derivation,  we use the fact
$$ 
\frac{t -r\sqrt{n}}{1-r^2}=\left(t -r\sqrt{n}\right)\cdot (1+O(r^2))=t -r\sqrt{n} +o(1)
$$ 
due to the assumptions $r \sqrt{n/\log p}\to 2-\sqrt{2}$
and $\log p = o(n^{1/5})$. From definition 
$\kappa_n:= (r \sqrt{n})(\log p)^{-1/2} -(2-\sqrt{2})$, we see $r \sqrt{n}=\kappa_n\sqrt{\log p}+(2-\sqrt{2})\sqrt{\log p}$. Use this identity and replace ``$t$'' in \eqref{whysha} with its definition in \eqref{toujiang} to obtain 
\begin{align*}
    t_1  &= \frac{x}{\sqrt{2}} - \sqrt{2} \kappa_n\log p + \lb \frac{1}{2} - \frac{1}{2 \sqrt{2}} \rb\log \log p + \lr \frac{1}{2}  - \frac{1}{2\sqrt{2}} \rr\log (4 \pi) +o(1)  \\
    & =- \lambda_n + \frac{x}{\sqrt{2}} + K_2+ o(1)
\end{align*}
where $\lambda_n=\sqrt{2} (\log p) \kappa_n + (8^{-1/2}-2^{-1})\log \log p$ and $K_2 = (2^{-1}-8^{-1/2})\log (4 \pi)$. Therefore, by \eqref{tata} and \eqref{mangong} we see that $\mb{P}(M^{(1)}_n \leq t)=\mb{P} (J_n \leq t_1) \to 1$ if $\lambda_n \rar -\infty$ and
$\mb{P}(M^{(1)}_n \leq t)$ goes to $\exp(-e^{-x/\sqrt{2}-K_2+\lambda})$ if 
$\lambda_n \rar \lambda$. Since $M_n=\max\{M^{(1)}_n, M^{(2)}_n\}$ and $M^{(1)}_n$ and $M^{(2)}_n$ are independent, it follows from \eqref{pianzi} that
\begin{eqnarray*}
\mb{P}(M_n \leq t) = \mb{P}(M^{(1)}_n \leq t)\cdot \mb{P}(M^{(2)}_n \leq t) \rar \exp(-e^{-x/\sqrt{2}-K_2+\lambda})\cdot e^{-K_1e^{-x}}
\end{eqnarray*}
if $\lambda_n \rar \lambda$. Also, it is easy to see $\mb{P}(M_n \leq t) \to \exp\{-K_1e^{-x}\}$ if $\lambda_n \rar -\infty$. These together with \eqref{earth2} imply statements  (i) and (ii) from the statement of Theorem \ref{large}. Certainly, if $\lambda_n \rar \infty$ then $t_1\to -\infty$ and hence $\mb{P}(M^{(1)}_n \leq t) \to 0$ by \eqref{mangong}. This implies $\mb{P}(M_n \leq t) \to 0$. In order to get a non-degenerate  limit, we need to change the scaling of $t$ in \eqref{toujiang}.

Now we treat the case $\lambda_n \rar \infty$. First, from the notation 
\begin{eqnarray*}
\kappa_n= \frac{r \sqrt{n}}{\sqrt{\log p}} -(2-\sqrt{2}) ~~~\mbox{and}~~~
\lambda_n=\sqrt{2} (\log p) \kappa_n + (8^{-1/2}-2^{-1})\log \log p,
\end{eqnarray*}
we get
\begin{eqnarray}\label{fuyu}
\lambda_n=r\sqrt{2n\log p}+(2-\sqrt{8})\log p +\left(\frac{1}{\sqrt{8}}-\frac{1}{2}\right)\log\log p.
\end{eqnarray}
Second, fix $x\in \mb{R}$, denote $u=u(x)$ by 
$$u= (1-r^2) \lr \frac{x}{ \sqrt{2\log p}} + \sqrt{2\log p} - \frac{\log \log p + \log (4 \pi)}{2 \sqrt{2\log p}} \rr + r\sqrt{n}.$$
From \eqref{tata} it is known $\mb{P}(M^{(1)}_n \leq u) \rar \exp\{-e^{-x}\}$. We shall prove $\mb{P}(M^{(2)}_n \leq u) \rar 1$. Write
\begin{align*}
    \mb{P}(M^{(2)}_n \leq u) &= \mb{P} \lr 2\sqrt{\log p} \lr M^{(2)}_n - 2\sqrt{\log p}+ \frac{\log \log p + \log (4 \pi)}{4 \sqrt{\log p}} \rr \leq u_1 \rr,
\end{align*}
where 
\begin{eqnarray*}
  u_1=   2\sqrt{\log p} \lr  u - 2\sqrt{\log p}+ \frac{\log \log p + \log (4 \pi)}{4 \sqrt{\log p}} \rr.
  \end{eqnarray*}
  Put the expression of $u$ into $u_1$ and use \eqref{fuyu} to see
  \begin{eqnarray*}
   u_1 &=& 2r\sqrt{n\log p}+(2\sqrt{2}-4)\log p +\left(\frac{1}{2}-\frac{1}{\sqrt{2}}\right)\log\log p+ O(1)\\
    &=& \sqrt{2}\lambda_n + O(1)\to \infty
\end{eqnarray*}
since $\lambda_n \rar \infty$. Thus $\mb{P}(M^{(2)}_n \leq u) \rar 1$ by \eqref{pianzi} and consequently
$$\mb{P}(M_n \leq u) =\mb{P} (M^{(1)}_n \leq u)\cdot \mb{P}(M^{(2)}_n \leq u) \rar e^{-e^{-x}}.$$ Finally, by \eqref{earth2},
$$
\mb{P} \lr \sqrt{n}\mc{L}_n+ O_{\mb{P}}\Big(\frac{\log p}{\sqrt{n}}\Big) \leq u \rr \rar e^{-e^{-x}}.
$$
Reorganize the event in the probability via the expression of $u$ and then apply the Slutsky lemma, we obtain that $c_n\mc{L}_n-d_n$ converges weakly to the  cdf $\exp\{-e^{-x}\}$, where
\begin{eqnarray*}
c_n = \frac{\sqrt{2n \log p}}{1-r^2}\ \ \ \mbox{and}\ \ \ d_n = \frac{r \sqrt{2n \log p}}{1-r^2} + 2\log p - \frac{1}{2} \big[ \log \log p + \log(4 \pi) \big]. 
\end{eqnarray*}
The proof is completed. \hfill$\square$

\subsection{Proof of Theorem \ref{Toeplitz}} In this section, the notation $H_{i,i+s}$ and $Z_{i,i+s}$ are used to indicate the coordinate of the random vectors $\bm{H}$ and $\bm{Z}$ at index $(i,i+s)$, respectively. If there is a sequence of such random vectors, we shall denote them by $H^{(n)}_{i,i+s}$ and $Z^{{ n}}_{i,i+s}$, respectively, with $n=1, 2, \cdots.$   We first collect a few useful lemmas whose proofs are presented in the supplementary material. 

\begin{lemma} \label{covariance}
Let $d\geq 1$ be a fixed integer and $\la r_k \ra_{k=1}^{\infty}$ be in \eqref{cov2} with $r_1=\cdots =r_d$. 
Let  $x_1, x_2, .., $ be an infinite sequence of random variables with $x_i \sim N(0, 1)$ for each $i$ and $\mbox{Cov}(x_i, x_j)=r_{|i-j|}$ for any $1\leq i< j$. Define
$$
Q_{i,i+s}  =\frac{r_1}{2(1-r_1^2)}(x_{i}^2-x_{i+s}^2)+ \frac{(x_{i+s} -r_1x_i)x_{i}}{1-r_1^2}
$$
for $i\geq 1$ and $1\leq s\leq d.$ Then $\mbox{Var}(Q_{i,i+s})=1$ and  $\mb{E}( Q_{i,i+1}Q_{j,j+1})$ is equal to
 $$
 \frac{1}{(1-r_1^2)^2} \Big[ r_1^2 r_k^2 + \frac{1}{2}r_1^2(r_{k-1}^2 + r_{k+1}^2) + r_k^2 + r_{k-1}r_{k+1} -2r_1r_k(r_{k-1}+r_{k+1}) \Big]$$
 for $k:=|i-j| \geq 1$ and $1\leq s, t\leq d$. 
 Furthermore, $\mb{E}( Q_{i,i+s}Q_{j,j+t})$  is a quadratic polynomial of $r_{|i-j+\alpha s+\beta t|}$ with $\alpha, \beta \in \{-1, 0, 1\}$.  Finally, 
if $\lim_{m\to\infty}r_{m} \sqrt{\log m} = \gamma \in [0, \infty)$, then 
$$\sup_{i,j,s,t}\Big| \mb{E}\left( Q_{i,i+s}Q_{j,j+t} \right) \log |i-j| - \frac{2\gamma^2}{(1+r_1)^2} \Big| \to 0$$ 
as $m \rightarrow \infty$, where the supremum is taken over $i\geq 1$, $j\geq 1$, $|i-j|\geq m$ and $1\leq s, t\leq d$. 
\end{lemma}

\begin{lemma} \label{<1}
 Assume the same setting as in Lemma \ref{covariance}. 
Then 
    $\sup_{i,j,s,t} \left| \mb{E}\left( Q_{i,i+s}Q_{j,j+t} \right) \right| < 1,$ 
    where the supremum is taken over $i\geq 1$, $j\geq 1$, $|i-j|\geq m$ and $1\leq s, t\leq d$. 

\end{lemma}

\begin{lemma}\label{mixed}
Let $\{\kappa_n \geq 1;\, n\geq 1\}$ be a sequence of integers with $\lim_{n\to\infty}\kappa_n=\infty$. 
For each $n\geq 1$ let $\{  X^{(n)}_i;\, 1 \leq i \leq \kappa_n\}$ be a  (possibly nonstationary) sequence of standard normals with covariance $r^{(n)}_{ij}=\mb{E}(X^{(n)}_iX^{(n)}_j)$.   Assume    $\sup_{i,j, n: 1\leq i<j \leq \kappa_n} |r^{(n)}_{ij}| < 1$. Set $M_{n}=\max_{1 \leq i \leq \kappa_n} X^{(n)}_i$.  If there exists $\gamma\geq 0$ such that
\begin{eqnarray*}
\sup_{1\leq i,j\leq \kappa_n, |i-j|\geq k} \left| r^{(n)}_{ij} \log |i-j| - \gamma \right| \to 0
\end{eqnarray*}
as $n\to \infty$ and $k\to\infty$ with $k<\kappa_n$, then
\begin{eqnarray*}
\sqrt{2\log \kappa_n} \left[M_n - \Big(\sqrt{2 \log \kappa_n} - \frac{\log \log \kappa_n + \log (4\pi)}{2\sqrt{2 \log \kappa_n}}\Big) \right]
\end{eqnarray*}
converges weakly to $-\gamma+\sqrt{2\gamma} Z + G$, where $Z \sim N(0,1)$, $G$ has cdf  $\exp(-e^{-x})$  and the two random variables are independent.
\end{lemma}

\begin{lemma}\label{niceyinli} 
Let $r_0=1, r_1, r_2, \cdots $ be non-negative constants. Assume the matrix $\bm{\Sigma}_k:=(r_{|i-j|})_{1\leq i, j\leq k}$ is non-negative  definite for each $k\geq 2$. Then there exists an infinite sequence of random variables $x_1, x_2, \cdots $ such that they are jointly normal, $x_i\sim N(0, 1)$ for each $i\geq 1$ and $(x_1, \cdots, x_k)^T$ has covariance matrix $\bm{\Sigma}_k$ for each $k\geq 1$.
\end{lemma}

\noindent\textbf{Proof of Theorem \ref{Toeplitz}}. Recall 
\eqref{language} that 
$$
\mc{L}_n =\max_{1\leq i<j \leq p} \hat{\rho}_{i,j}.
$$
The argument in the proof is similar to that of the proof of  Theorem \ref{slowdecay}. We will only present the necessary changes. 

{\it Step 1: Reduction of $\mc{L}_n$ to a maximum over a smaller subset.}    One can show that, with probability going to $1$, we have
    $$\max_{\substack{ 1 \leq i < j \leq p\\ {1 \leq j-i \leq d}}} \hat{p}_{i,j} \geq \max_{ \substack{ 1\leq i<j \leq p\\ j - i \geq d+1}} \hat{p}_{i,j}$$
    as $n\to\infty$.
    The technical details are almost identical to that of {\it Step 1} in the proof of Theorem \ref{slowdecay} and we skip it for brevity.
    
    {\it Step 2: Linearize $\hat{\rho}_{i,j}$ in \eqref{wangshi}}. Define $I=\la (i,i+s):1\leq i <i+s \leq p, 1\leq s \leq d \ra$.  By the same argument as in the {\it Step 2} of the proof of Theorem \ref{slowdecay}, we have
    \begin{eqnarray*}
       && \max_{(i,i+s) \in I} \sqrt{n} \hat{p}_{i,i+s}\\
       &=& \max_{(i,i+s) \in I} \lb \frac{1}{\sqrt{n}} \sum_{k=1}^{n} x_{k,i}x_{k,i+s}\rb\cdot \lb 2-\frac{1}{2n}\sum_{k=1}^{n} (x_{k,i}^2+x_{k,i+s}^2)\rb +O_{\mb{P}}\Big(\frac{\log p}{\sqrt{n}}\Big).
    \end{eqnarray*}
    Define 
    \begin{eqnarray*}
    H^{n}_{i,i+s} = \lb\frac{1}{\sqrt{n}} \sum_{k=1}^{n} x_{k,i}x_{k,i+s}\rb\cdot \lb 2-\frac{1}{2n}\sum_{k=1}^{n} (x_{k,i}^2+x_{k,i+s}^2)\rb.
    \end{eqnarray*}
    Then $\max_{(i,i+s) \in I} \sqrt{n} \hat{p}_{i,i+s}=\max_{(i,i+s) \in I}  H^{n}_{i,i+s}.$ 
    Note that the index set $I$ in Lemma \ref{covariance} has asymptotically $pd$ elements and hence the discussions of {\it Step 1} of the proof of Theorem \ref{slowdecay} can be adapted with minor changes. The details are also skipped. In summary, from {\it Step 1} and {\it Step 2} above, we have 
    \begin{eqnarray}\label{dashun}
    \mc{L}_n =\max_{(i,i+s) \in I}  H^{n}_{i,i+s} \to 1~~ \mbox{in probability}.
    \end{eqnarray}
    
    {\it Step 3: A decomposition of $H^{n}_{i,i+s}$}. 
    Recall $\textbf{x}_1, \textbf{x}_2, \cdots, \textbf{x}_n$ is a random sample from $N(\bm{0}, \bm{\Sigma})$ with $\bm{\Sigma}$ in \eqref{cov1} and $\textbf{x}_i=(x_{i1},x_{i2},\cdots,x_{ip})^T$ for each $i$. Hence the rows of data matrix $(\textbf{x}_1, \textbf{x}_2, \cdots, \textbf{x}_n)^T=(x_{ij})_{n\times p}$ are i.i.d. random vectors with distribution $N(\bm{0}, \bm{\Sigma})$.
    For $1\leq s \leq d$ and $(i,i+s) \in I$, we denote by $\xi^{k}_{i,s}$ the Gaussian random variable such that
    $$\xi^{s}_{k,i}=x_{k,i+s}-r_1 x_{k,i}.$$
    It's easy to see  $\xi_{k, i}^{s} \sim N(0,1-r_1^2)$ and $\xi_{k, i}^{s}$ is independent of $x_{k,i}$.     Fix $(i,i+s) \in I$, we have 
  \begin{align*} 
    H^{n}_{i,s} &= \lb\frac{1}{\sqrt{n}} \sum_{k=1}^{n} x_{k,i}x_{k,i+s}\rb \lb 2-\frac{1}{2n}\sum_{k=1}^{n} (x_{k,i}^2+x_{k,i+s}^2)\rb \\
    &= \left[\frac{r_1}{\sqrt{n}}\sum_{k=1}^{n}(x_{k,i}^2-1)+\frac{1}{\sqrt{n}}\sum_{k=1}^{n}\xi_{k, i}^{s}x_{k,i}+r_1\sqrt{n}\right] \left[1-\frac{1}{2n}\sum_{k=1}^{n} (x_{k,i}^2+x_{k,i+s}^2-2)\right]\\ \nonumber
    & = \frac{r_1}{\sqrt{n}}\sum_{k=1}^{n}(x_{k,i}^2-1) - \frac{r_1}{2n^{3/2}}\sum_{k=1}^{n} (x_{k,i}^2+x_{k,i+s}^2-2)\cdot\sum_{k=1}^{n}(x_{k,i}^2-1)\\ \nonumber
    &+ \frac{1}{\sqrt{n}}\sum_{k=1}^{n}\xi_{k, i}^{s}x_{k,i} - \frac{1}{2n^{3/2}}\sum_{k=1}^{n} (x_{k,i}^2+x_{k,i+s}^2-2)\cdot\sum_{k=1}^{n}\xi_{k, i}^{s}x_{k,i}\\ \nonumber
    &+ r_1\sqrt{n}-\frac{r_1}{2\sqrt{n}}\sum_{k=1}^{n} (x_{k,i}^2+x_{k,i+s}^2-2)\\ \nonumber
    &:=A_n+B_n+C_n+D_n+r_1\sqrt{n}+E_n.
\end{align*}
    With almost identical arguments to Step 3 of Theorem \ref{slowdecay}, one can show that $B_n,D_n$ are of order $O_{\mb{P}}(n^{-1/2}\log p)$. Therefore, by combining $A_n$,  $C_n$ and $E_n$ together, we obtain
    $$\max_{(i,i+s) \in I} H^{n}_{i,i+s} -r_1 \sqrt{n} = \max_{(i,i+s) \in I} \frac{1}{\sqrt{n}} \sum_{k=1}^{n} \lb \frac{r_1}{2}(x_{ki}^2-x_{k,i+s}^2)+\xi_{k, i}^{s}x_{k,i} \rb + O_{\mb{P}}\Big(\frac{\log p}{\sqrt{n}}\Big).$$
    Define 
    \begin{eqnarray}\label{henhaoba}
     M_n=\max_{(i,i+s) \in I} \frac{1}{\sqrt{n}} \sum_{k=1}^{n} \Big[\frac{r_1}{2}\big(x_{ki}^2-x_{k,i+s}^2)\big)+\xi^{s}_{k,i}x_{k,i} \Big].
    \end{eqnarray}
    By \eqref{dashun}, we have 
    \begin{eqnarray}\label{shengli}
    \sqrt{n}\mc{L}_n -r_1 \sqrt{n} =M_n + O_{\mb{P}}\Big(\frac{\log p}{\sqrt{n}}\Big)
    \end{eqnarray}
    with probability going to $1$ as $n\to\infty$.

     {\it Step 4: Gaussian approximation.} We now examine $M_n$ in  \eqref{henhaoba}. 
    Define a sequence of i.i.d. $|I|$-dimensional random vectors $\la \bm{Q}_k \ra_{k=1}^{n}$ such that 
      $$ \bm{Q}_k= \frac{1}{1-r_1^2}\left(\frac{r_1}{2}\big(x_{k,i}^2-x_{k,i+s}^2\big)+\xi^{s}_{k,i}x_{k,i} \right)_{(i,i+s) \in I}.$$
    Let $\bm{S}^{(n)}=n^{-1/2}\sum_{k=1}^{n} \bm{Q}_k$. Note that we have the relation
    \begin{eqnarray}\label{saoxue}
      \frac{M_n}{1-r_1^2}=\max_{(i,i+s)\in I} S^{(n)}_{i,i+s},
      \end{eqnarray}
      where, as usual, $S^{(n)}_{i,i+s}$ denotes the $(i, i+s)$-coordinate of the vector $\bm{S}^{(n)}$. Notice $\bm{Q}_1 \cdots, \bm{Q}_n$ are i.i.d. random vectors with mean $\bm{0}$. Reviewing the notation $Q_{i,i+s}$ in Lemma \ref{covariance}, we see
      \begin{eqnarray}\label{Nash}
      \mbox{Var}\left(\bm{Q}_1(i, i+s)\right)=1~~ \mbox{and}~~\mbox{Cov}\left(\bm{Q}_1(i, i+s), \bm{Q}_1(j, j+t)\right)= \mb{E}( Q_{i,i+s}Q_{j,j+t})
      \end{eqnarray}
   for $i\ne j$.  By Lemma \ref{vershynin}, we have
     \begin{align*}
         \sup_{(i,i+s) \in I} \left\|  \frac{r_1}{2}\big(x_{k,i}^2-x_{k,i+s}^2\big)+\xi^{s}_{k,i}x_{k,i} \right\|_{\psi_1} & \leq C\cdot \Big( \| N(0,1) \|_{\psi_2}^2 + \sup_{i,k,s} \| \xi^{s}_{k,i} \|_{\psi_2} \cdot \| N(0,1) \|_{\psi_2} \Big) \\
         & \leq C\cdot\| N(0,1) \|_{\psi_2}^2 
     \end{align*}
     where $C$ is the constant in Lemma \ref{vershynin} and we use the fact that $\xi^{s}_{k,i} \sim N(0,1-r_1^2)$. Let 
     $\bm{Z}^{(n)}=(Z^{(n)}_{i,i+s})_{(i,i+s)\in I}$ be a $|I|$-dimensional Gaussian vector with the same covariance structure as in \eqref{Nash}, that is, 
    \begin{eqnarray}\label{Nash_da}
      \mbox{Var}\left(Z^{(n)}_{i,i+s}\right)=1~~ \mbox{and}~~\mbox{Cov}\left(Z^{(n)}_{i,i+s}, Z^{(n)}_{j, j+t}\right)= \mb{E}( Q_{i,i+s}Q_{j,j+t}).
      \end{eqnarray}
     Apply the high-dimensional CLT (see Theorem S.2 from  \cite{JiangPham21} and also \cite{Koike}) to obtain that
    \begin{eqnarray*}
    \sup_{t \in \mb{R}^{|I|}} \bigm| \mb{P}(\bm{S}^{(n)} \leq t) - \mb{P}(\bm{Z}^{(n)} \leq t) \bigm| \leq O\left(\frac{(\log |I|)^{5/6}}{n^{1/6}}\right)=O\left(\frac{(\log p)^{5/6}}{n^{1/6}}\right),
    \end{eqnarray*}
    where the last identity  holds since $|I| = O(p)$ by the definition of $I$. In particular, from the above identity together with \eqref{shengli} and  \eqref{saoxue} we conclude
    \begin{eqnarray}\label{zaosheng}
    \sup_{s \in \mb{R}} \left| \mb{P}\left(\frac{\sqrt{n} \mc{L}_n -r_1 \sqrt{n}}{1-r_1^2}+ O_{\mb{P}}\left(\frac{\log p}{\sqrt{n}}\right)\leq s\right) - \mb{P}(\mathcal{E}_n \leq s) \right| =O\left(\frac{(\log p)^{5/6}}{n^{1/6}}\right)
    \end{eqnarray}
    where $\mathcal{E}_n:=\max_{(i,i+s)\in I} Z^{n}_{i,i+s}$. We will derive  the asymptotic distribution of  $\mathcal{E}_n$ next. 
    
    {\it Step 5: The asymptotic distribution of $\mathcal{E}_n$.} 
    We will apply Lemma \ref{mixed} to deduce the result. To this end, write the coordinates of $\bm{Z}^{(n)}$  (which is an unbalanced triangular array) by
\begin{eqnarray}
\begin{array}{ c c c c c c c c }
Z^n_{1,2} & Z^n_{2,3} & Z^n_{3,4} & \cdots  & Z^n_{p-d,p-d+1} & \cdots  & Z^n_{p-2,p-1} & Z^n_{p-1,p}\\
Z^n_{1,3} & Z^n_{2,4} & Z^n_{3,5} & \cdots  & Z^n_{p-d,p-d+2} & \cdots  & Z^n_{p-2,p} & \\
\vdots  & \vdots  & \vdots  &  & \vdots  &  &  & \label{jieqilai}\\
Z^n_{1,d+1} & Z^n_{2,d+2} & Z^n_{3,d+3} & \cdots  & Z^n_{p-d,p}, &  &  & 
\end{array}
\end{eqnarray}
where we write  $Z^{n}_{i,i+s}$ for $Z^{(n)}_{i,i+s}$ in order to simplify  notation. Observe the triangular array has two parts: a matrix of $d$ rows and $p-d$ columns in the left side and a triangle in the right side with  total number of entries $1+2+\cdots +(d-1)=(1/2)d(d-1)$. Denote $m_n=\max_{1\leq i \leq p-d,\, 1\leq s \leq d}Z^{n}_{i,i+s}$. Since $Z^{n}_{i,i+s} \sim N(0, 1)$, we have
\begin{eqnarray*}
\mb{P}(m_n> a) \leq \mb{P}(\mathcal{E}_n > a)\leq \mb{P}(m_n > a)+\frac{1}{2}d(d-1)\cdot\mb{P}(N(0, 1)>a)
\end{eqnarray*}
for any $a \in \mathbb{R}.$ In the following discussion we will see $a=a_n\to\infty$, then the last probability above goes to zero since $d$ is free of $n$. Hence, for any sequence of $a_n>0$ with $\lim_{n\to\infty}a_n=\infty$ we have $\lim_{n\to\infty}\left[\mb{P}(\mathcal{E}_n > a_n)-\mb{P}(m_n > a_n) \right]=0$, or equivalently, 
$\lim_{n\to\infty}\left[\mb{P}(\mathcal{E}_n \leq a_n)-\mb{P}(m_n\leq a_n) \right]=0.$
We claim that, as $n\to \infty$, 
\begin{eqnarray} 
&&    \mb{P}\left( \sqrt{2 \log (pd)} \left(m_n- \sqrt{2 \log (pd)}+\frac{\log \log (pd)}{2\sqrt{2 \log (pd)}} \right)\leq x \right) \nonumber\\
&\rightarrow & \mb{P}\Big( -\gamma_0+\sqrt{ 2\gamma_0}Z+G \leq x + \frac{1}{2}\log (4\pi)\Big),\label{fact9} 
\end{eqnarray}
where $\gamma_0=2\gamma^2(1+r_1)^{-2},
~~Z \sim N(0,1),~~ G \ \mbox{has}\ \mbox{cdf}\  \exp(-{e^{-x}})$,   $Z$ and $G$ are independent. 
Assuming \eqref{fact9} holds, then it also holds if ``$m_n$'' is replaced with ``$\mathcal{E}_n$''. This implies
$$
\sqrt{2 \log (pd)} \left(\mathcal{E}_n- \sqrt{2 \log (pd)}+\frac{\log \log (pd)}{2\sqrt{2 \log (pd)}} \right)
$$
converges weakly to $-\gamma_0+\sqrt{2\gamma_0}Z+G -(1/2)\log (4\pi)$. 
This, \eqref{zaosheng}, the assumption $\log p = o(n^{1/5})$ and the Slutsky lemma yield 
\begin{eqnarray*}
\sqrt{2 \log (pd)} \left(\frac{\sqrt{n} \mc{L}_n -r_1 \sqrt{n}}{1-r_1^2}- \sqrt{2 \log (pd)}+\frac{\log \log (pd)}{2\sqrt{2 \log (pd)}}\right) 
\end{eqnarray*}
converges weakly to 
$-\gamma_0 + \sqrt{2 \gamma_0} Z+G -(1/2)\log (4\pi)$. This gives part (i). Now we prove \eqref{fact9} in the following. 

Recall notation $Z^{n}_{i,i+s}=Z^{(n)}_{i,i+s}$. We know  $(Z^{n}_{i,i+s})_{(i,i+s)\in I}$ is a $|I|$-dimensional Gaussian vector with mean zero and covariance structure presented in \eqref{Nash_da}. Recollecting $Q_{i,i+s}$ in Lemma \ref{covariance}, we see from the lemma that $\mb{E}(Z^{n}_{i,i+s}Z^{n}_{j,j+t})$ is free of $n$ and $p_n$, and 
\begin{eqnarray}\label{xingxingdd} 
\sup_{i,j,s,t}\Big| \mb{E}\left(Z^{n}_{i,i+s}Z^{n}_{j,j+t} \right) \log |i-j| - \frac{2\gamma^2}{(1+r_1)^2} \Big| \to 0
\end{eqnarray}
as $m \rightarrow \infty$,  where the supremum is taken over $i\geq 1$, $j\geq 1$, $|i-j|\geq m$ and $1\leq s, t\leq d$. We will apply Lemma \ref{mixed} to show \eqref{fact9}. To do so, we need to  transform 
the matrix  $\{Z^{n}_{i,i+s};\, 1\leq i \leq p-d,\, 1\leq s \leq d\}$ in \eqref{jieqilai} into a column. We now label the random variables as a sequence by assigning the order vertically, that is, we place the second column from the left in \eqref{jieqilai} at the bottom of the first column, and put the third column in \eqref{jieqilai} at the bottom of the previous vector, and so on. In this way, we obtain a $(p-d)d$-dimensional vector. Precisely, it is a map $\sigma$: $\{(i,i+s);\, 1\leq i \leq p-d,\, 1\leq s \leq d\} \to \{1,2, \cdots, (p-d)d\}$ satisfying $\sigma(i,i+s)=(i-1)d+s$. Observe that $|\mb{E}( Z^n_{i,i+s}Z^n_{j,j+t})|\leq 1$ and 
\begin{eqnarray}\label{diwen}
\log |\sigma(i,i+s)- \sigma(j,j+t)|
&=&\log |(i-j)d+(s-t)|\nonumber\\
& = & \log |i-j| + \log \left|d+\frac{s-t}{i-j}\right|.
\end{eqnarray}
Denote $\rho(s,t,i,j)= \log |d+(s-t)(i-j)^{-1}|$. Then $|\sup\rho(s,t,i,j)| <\infty$ where the supremum is taken over $i, j, s, t$ with $|i-j|>2d$, $1\leq s, t \leq d.$ The assertion  \eqref{xingxingdd} implies that 
\begin{eqnarray}\label{rezhene}
\sup_{i,j,s,t}\Big| \mb{E}\left(Z^{n}_{i,i+s}Z^{n}_{j,j+t} \right)\Big| \to 0
\end{eqnarray}
as $m \rightarrow \infty$,  where the supremum is taken over $i\geq 1$, $j\geq 1$, $|i-j|\geq m$ and $1\leq s, t\leq d$. This joining \eqref{xingxingdd}, \eqref{diwen} and \eqref{rezhene} implies 
\begin{align*} 
  \sup_{i,j,s,t}  \Big| \mb{E}\left( Z^n_{i,i+s}Z^n_{j,j+t} \right) \log |\sigma(i,i+s)-\sigma(j,j+t)| - \frac{2\gamma^2}{(1+r_1)^2} \Big|=0
\end{align*}
as $m\to\infty$, where the supremum runs over $i,j,s,t$ satisfying $1\leq i, j \leq p-d$, $1\leq s, t \leq d$ and $|i-j|\geq m$. Recall $p=p_n$ and $m_n=\max_{1\leq i \leq p-d,\, 1\leq s \leq d}Z^{n}_{i,i+s}$. For each $n\geq 2$, take $\kappa_n=(p-d)d$ in Lemma \ref{mixed} and use \eqref{Nash_da} and Lemma \ref{<1}  to see
\begin{eqnarray*}
\sqrt{2\log \kappa_n} \left[m_n - \Big(\sqrt{2 \log \kappa_n} - \frac{\log \log \kappa_n + \log (4\pi)}{2\sqrt{2 \log \kappa_n}}\Big) \right]
\end{eqnarray*}
converges weakly to $-\gamma_0+\sqrt{2\gamma_0} Z + G$, where $\gamma_0=2\gamma^2(1+r_1)^{-2}$,  $Z \sim N(0,1)$, $G$ has cdf  $\exp(-e^{-x})$  and the two random variables are independent. Trivially $\log \kappa_n=\log (pd) + O(p^{-1})$. Then the above weak convergence is still true if ``$\kappa_n$'' is replaced with ``$pd$'' via the Slutsky lemma. This implies \eqref{fact9}.  

(ii) Review \eqref{zaosheng}. Similar to the proof of (i), we only need to analyze $\mathcal{E}_n$ when $d=1$. In this case, $\mathcal{E}_n =  \max_{1 \leq i \leq p-1} Z^{n}_{i,i+1}$. Lemma \ref{covariance} and \eqref{Nash_da} imply that $\{ Z^{n}_{i,i+1} \}_{1 \leq i \leq p-1}$ is a stationary sequence of $N(0, 1)$-distributed random variables with correlation sequence $\{f(k);\, 1\leq k\leq p-2\}$ given by $f(0)=1$ and 
$$f(k)= \frac{1}{(1-r_1^2)^2} \Big[ r_1^2 r_k^2 + \frac{1}{2}r_1^2(r_{k-1}^2 + r_{k+1}^2) + r_k^2 + r_{k-1}r_{k+1} -2r_1r_k(r_{k-1}+r_{k+1}) \Big]$$
for all $k \geq 1$. In particular, $(f(|i-j|))_{1\leq i, j\leq k}$ is non-negative definite for each $k\geq 2$. Considering the infinite sequence $\{f(0), f(1), f(2), \cdots \}$, by Lemma \ref{niceyinli}, there exists an infinite sequence of random variables $\xi_1, \xi_2, \cdots $ such that they are jointly normal, $\xi_i\sim N(0, 1)$ for each $i\geq 1$ and $(\xi_1, \cdots, \xi_k)^T$ has covariance matrix $(f(|i-j|))_{1\leq i, j\leq k}$ for each $k\geq 1$.  By assumption, $\{f(k);\, k\geq 1\}$ is monotone decreasing for large $k$,   $\lim_{k\to\infty}f(k)= 0$ and $f(k) \log k$ increases to $\infty$. From Theorem 6.6.4 in \cite{Leadbetter}  we deduce that
$$
\frac{1}{\sqrt{f(k)}}\left[\max\{\xi_1, \cdots, \xi_k\}-\sqrt{1-f(k)}\,b_k\right] \to N(0, 1), 
$$
 where 
$
b_k=\sqrt{2\log k}-\frac{\log \log k+\log (4\pi)}{2\sqrt{2\log k}}.
$
Obviously, $(Z^{n}_{1,2}, \cdots, Z^{n}_{p-1, p})$ and $(\xi_1, \cdots, \xi_{p-1})$ have the same distribution. It follows that
\begin{eqnarray*}
\frac{1}{\sqrt{f(p-1)}}\left[\mathcal{E}_n - b_{p-1}\sqrt{1-f(p-1)}\right]
\end{eqnarray*}
converges weakly to $N(0, 1)$. From \eqref{zaosheng} we see that
\begin{eqnarray*}
\frac{1}{\sqrt{f(p-1)}}\left[\frac{\sqrt{n} \mc{L}_n -r_1 \sqrt{n}}{1-r_1^2}+ O_{\mb{P}}\left(\frac{\log p}{\sqrt{n}}\right) - b_{p-1}\sqrt{1-f(p-1)}\right] \to N(0, 1)
\end{eqnarray*}
weakly. It is checked in Section 3.6 from the supplement Jiang and Pham \cite{JiangPham21} that
\begin{eqnarray*}
\frac{1}{\sqrt{f(p-1)}}\left( O_{\mb{P}}\left(\frac{\log p}{\sqrt{n}}\right) +(b_p- b_{p-1})\sqrt{1-f(p-1)}\right) \to 0
\end{eqnarray*}
as $n\to \infty$. The two assertions above and the Slutsky lemma conclude the proof. \hfill$\square$

\medskip

\noindent\textbf{Supplement to “Asymptotic Distributions of Largest Pearson Correlation Coefficients under Dependent Structures”} contains proofs and further discussions.



\begin{thebibliography}{}

\bibitem{Bentkus}
BENTKUS, V. (2003). On the dependence of the Berry–Esseen bound on dimension. \textit{Journal of Statistical Planning and Inference 113,   385–402.}

\bibitem{Cai17}
CAI, T. (2017). Global testing and large-scale multiple testing for high-dimensional covariance structures, \textit{Annual Review of Statistics and Its Application  Volume 4, 423-446.}


\bibitem{Jiang13}
CAI, T., FAN, J. and JIANG, T.  (2013). Distributions of angles in random packing on spheres, \textit{J. Mach. Learn. Res. 14 1837–1864.}

\bibitem{CaiMa13}
CAI, T. and MA, Z. (2013). Optimal hypothesis testing for high dimensional covariance matrices, \textit{Bernoulli 19(5B): 2359-2388.}

\bibitem{Jiang12}
CAI, T. and JIANG, T. (2012). Limiting Laws of Coherence of Random Matrices with Applications to Testing Covariance Structure and Construction of Compressed Sensing Matrices, \textit{Ann. Stat. 39(3), 1496-1525.}

\bibitem{JiangCai12}
 CAI, T. and JIANG, T. (2012). Phase transition in limiting distributions of coherence of high-dimensional random matrices, \textit{J. Multivariate Anal. 107 24–39.}

\bibitem{Cai2011}
 CAI, T. and LIU, W. (2011). Adaptive thresholding for sparse covariance matrix estimation, \textit{Journal of the American Statistical Association 106 672–684.}

\bibitem{Cai2013}
CAI, T., LIU, W. and XIA, Y. (2013). Two-sample covariance matrix testing and support recovery in high-dimensional and sparse settings, \textit{Journal of the American Statistical Association 108 265–277.}

\bibitem{Cai2014}
CAI, T., LIU, W. and XIA, Y. (2014). Two-sample test of high dimensional means under dependence, \textit{Journal of the Royal Statistical Society Series B-statistical Methodology 76 349–372.}

\bibitem{Cai16}
CAI, T. and ZHANG, A.  (2016). Inference for high-dimensional differential correlation matrices, \textit{Journal of Multivariate Analysis 143 107–126.}

\bibitem{Kato20}
CHERNOZHUKOV, V.,  CHETVERIKOV, D. and KOIKE, Y. (2020). Nearly optimal central limit theorem and bootstrap approximations in high dimensions, \textit{arXiv:2012.09513.}

\bibitem{Kato19}
CHERNOZHUKOV, V.,  CHETVERIKOV, D.,  KATO, K. and KOIKE, Y. (2022). Improved central limit theorem and bootstrap approximations in high dimensions, \textit{Ann. Statist. 50(5), 2562-2586.}

\bibitem{Kato17}
CHERNOZHUKOV, V.,  CHETVERIKOV, D. and  KATO, K. (2017). Central limit theorems and bootstrap in high dimensions, \textit{Ann. Probab. 45, 2309-2352.}

\bibitem{kato14-anti}
CHERNOZHUKOV, V.,  CHETVERIKOV, D. and  KATO, K. (2015). Comparison and anti-concentration bounds for maxima of Gaussian random vectors, \textit{Probability Theory and Related Fields 162 (1-2), 47–70.}

\bibitem{Kato13}
CHERNOZHUKOV, V.,  CHETVERIKOV, D. and  KATO, K. (2013). Gaussian approximations and multiplier bootstrap for maxima of sums of high-dimensional random vectors, \textit{Ann. Statist. 41, 2786-2819.}

\bibitem{ChenLiu}
CHEN, L. and LIU, W. (2017). Testing independence with high-dimensional correlated samples, \textit{Ann. Statist. 46(2), 866-894.}
\bibitem{Jiang19}
FAN, J. and JIANG, T. (2019). Largest entries of sample correlation matrices from equi-correlated normal populations, \textit{Ann. Probab. 47(5): 3321-3374.}

\bibitem{FengJiang}
FENG, L.,  JIANG, T., LIU, B. and XIONG, W. (2021). Max-sum tests for cross-sectional dependence of high-demensional panel data, \textit{Ann. Stat. 50(2), 1124-1143.}

\bibitem{Gotze}
GOTZE, F. (1991). On the rate of convergence in the multivariate CLT, \textit{Ann. Probab. 19(2): 724-739.}


\bibitem{Jiang04}
JIANG, T. (2004). The asymptotic distributions of the largest entries of sample correlation matrices, \textit{Ann. Appl. Probab,14(2): 865-880.}

\bibitem{Jiang2019}
JIANG, T. (2019). Determinant of sample correlation matrix with application, \textit{Ann. Appl. Probab, 29(3): 1356-1397.}

\bibitem{JiangPham21}
JIANG, T. and PHAM, T. Supplement to “Asymptotic distributions of largest Pearson correlation coefficients under dependent structures”.

\bibitem{Koike}
KOIKE, Y. (2020). Notes on the dimension dependence in high-dimensional central limit theorems for hyperrectangles, \textit{Japanese Journal of Statistics and Data Science.}

\bibitem{D.Li}
LI, D. and XUE, L. (2015). Joint limiting laws for high-dimensional independence tests, \textit{arXiv:1512.08819.}

\bibitem{D.Li2}
LI, D. and XUE, L. and ZOU, H. (2018)
Applications of Peter Hall's Martingale Limit Theory to Estimating and Testing High Dimensional Covariance Matrices. 
\textit{Statistica Sinica, 28: 2657-2670.}

\bibitem{Li1}
LI, D., LIU, W. and ROSALKY, A. (2010). Necessary and sufficient conditions for the asymptotic distribution of the largest entry of a sample correlation matrix, \textit{Probability Theory and Related Fields 148, 5–35.}

\bibitem{Li2}
LI, D., QI, Y. and ROSALKY, A. (2012). On Jiang’s asymptotic distribution of the largest entry of a sample correlation matrix, \textit{J. Multivariate Anal. 111 256–270.}

\bibitem{Li3}
LI, D. and ROSALKY, A. (2006). Some strong limit theorems for the largest entries of sample correlation matrices, \textit{Ann. Appl. Probab. 16 423–447.}

\bibitem{Leadbetter}
LINDGREN, G., ROOTZÉN, H. and LEADBETTER, M. R. (1983). Extremes and related properties of random sequences and processes, \textit{Springer Series in Statistics.}

\bibitem{Liu16}
LIU, W., LIN, Z. and SHAO, Q. M. (2008). The asymptotic distribution and Berry–Esseen bound of a new test for independence in high dimension with an application to stochastic optimization, \textit{Ann. Appl. Probab. 18(6): 2337-2366.}



\bibitem{Portnoy}
PORTNOY, S. (1986). On the central limit theorem in $\mb{R}^p$ when $p \rightarrow \infty$, \textit{Probability Theory and Related Fields 73, 571–583.}

\bibitem{Polya}
SASVÁRI, Z. (1998). On a classical theorem in the theory of Fourier Integrals, \textit{Proc. Amer. Math. Soc. 126(3), 711–713.}

\bibitem{ShaoZhou}
SHAO, Q. M. and ZHOU, W. X. (2014). Necessary and sufficient conditions for the asymptotic distributions of coherence of ultra-high dimensional random matrices, \textit{Ann. Probab. 42 623–648.}


\bibitem{MW2019}
WRAINWRIGHT, J. M. (2019). \textit{High-dimensional statistics: A non-asymptotic viewpoint}. Cambridge University Press. 


\bibitem{Yu1}
Yu, X., Li, D. and Xue, L. (2022). Fisher's combined probability test for high-dimensional covariance matrices, \textit{Journal of the American Statistical Association, to appear.}

\bibitem{Yu2}
Yu, X., Li, D., Xue, L. and Li, R. (2022). Power-enhanced simultaneous test of high-dimensional mean vectors and covariance matrices with application to gene-set testing, \textit{Journal of the American Statistical Association, to appear.}

\bibitem{Zhou07}
ZHOU, W. (2007). Asymptotic distribution of the largest off-diagonal entry of correlation matrices, \textit{Trans. Amer. Math. Soc. 359 5345–5363.}

\end{thebibliography}
\end{document}